\documentclass[11pt,amsfonts]{amsart}

\usepackage{fullpage}
\usepackage{amsmath}
\usepackage{amsfonts}
\usepackage{latexsym}
\usepackage{graphicx}

\def\endproof{\qed \medskip}
\def\blacksquare{\hbox to .60em {\vrule width .60em height .60em}}

\newtheorem{theorem}{Theorem}[section]

\newtheorem{corollary}[theorem]{Corollary}

\newtheorem{lemma}[theorem]{Lemma}

\newtheorem{proposition}[theorem]{Proposition}

\newtheorem{remark}[theorem]{Remark}

\begin{document}

\title[]{On the Bartnik extension problem for the static vacuum Einstein equations}

\author[]{Michael T. Anderson and Marcus A. Khuri}

\address{Dept.~of Mathematics, Stony Brook University, Stony Brook, NY 11790}
\email{anderson@math.sunysb.edu}
\email{khuri@math.sunysb.edu}

\thanks{The first author is partially supported by NSF Grant DMS 0905159 and 1205947 while the
second author is partially supported by NSF Grant DMS 1007156 and a Sloan Research
Fellowship. \\
PACS 2010: 04.20.-q, 04.20.Cv, 02.40.Vh, 02.30.Jr}

\abstract{We develop a framework for understanding the existence of asymptotically flat
solutions to the static vacuum Einstein equations on $M = {\mathbb R}^{3}\setminus B$
with geometric boundary conditions on $\partial M \simeq S^{2}$. A partial existence
result is obtained, giving a partial resolution of a conjecture of Bartnik
on such static vacuum extensions. The existence and uniqueness of such
extensions is closely related to Bartnik's definition of quasi-local mass. }
\endabstract

\maketitle

\section{Introduction}

\setcounter{section}{1}
\setcounter{equation}{0}

  This paper is concerned with a conjecture of R. Bartnik \cite{B3}, \cite{B4} on the
existence and uniqueness of static solutions to the vacuum Einstein equations with
certain prescribed boundary data. On the physical side, this is closely related to the
issue of local mass in general relativity while, on the mathematical side, to the issue of
global existence and uniqueness for a rather complicated geometric non-linear system
of elliptic boundary value problems.

    Let $M$ be a 3-manifold diffeomorphic to ${\mathbb R}^{3}\setminus B$ where $B$ is
a 3-ball, so that $\partial M \simeq S^{2}$. The static vacuum Einstein equations are the
equations for a pair $(g, u)$ consisting of a smooth Riemannian metric $g$ on $M$
and a positive potential function $u: M \rightarrow {\mathbb R}^{+}$ given by
\begin{equation}\label{1.1}
uRic_{g} = D^{2}u, \ \ \Delta u = 0,
\end{equation}
where the Hessian $D^{2}$ and Laplacian $\Delta = tr D^{2}$ are taken with respect
to $g$. The equations \eqref{1.1} are equivalent to the statement that the 4-dimensional
metric
\begin{equation}\label{1.2}
g_{\mathcal{M}} = \pm u^{2}dt^{2} + g,
\end{equation}
on the 4-manifold $\mathcal{M} = {\mathbb R}\times M$ is Ricci-flat, i.e.
\begin{equation}\label{1.3}
Ric_{g_{\mathcal{M}}} = 0.
\end{equation}
This holds for either choice of sign in \eqref{1.2} and since most of the analysis of
the paper concerns the Riemannian data $(g, u)$ in \eqref{1.1}, we will assume $g_{\mathcal{M}}$
is Riemannian, and moreover identify $t$ in \eqref{1.2} periodically, to obtain
a metric on $\mathcal{M} = S^{1}\times M$ with $t$ replaced by the angular variable $\theta$.

  Given $(M, g, u)$ as above, let $\gamma$ be the Riemannian metric induced on
$S^{2} = \partial M$ and let $H$ be the mean curvature of $\partial M \subset (M, g)$,
(with respect to the inward unit normal into $M$). Then (one version of) the Bartnik 
conjecture \cite{B4} states that, given an arbitrary such pair in $C^{\infty}$,
\begin{equation}\label{1.4}
(\gamma, H) \in Met^{\infty}(S^{2})\times C_{+}^{\infty}(S^{2}), \ \ H > 0,
\end{equation}
there exists a unique asymptotically flat solution $(g, u)$ to the static vacuum
Einstein equations \eqref{1.1} inducing the boundary data $(\gamma, H)$ on $\partial M$.

  This conjecture is a natural outgrowth of Bartnik's concept of quasi-local mass
$m_{B}(\Omega)$, \cite{B2}, \cite{B3}, defined as follows. Let $(\Omega, g)$ be a
smooth compact 3-manifold with smooth boundary of non-negative scalar curvature, and
define an {\it admissible extension} of $(\Omega, g)$ to be a complete, asymptotically
flat 3-manifold $(\widetilde M, g)$ of non-negative scalar curvature in which
$(\Omega, g)$ embeds isometrically and is not enclosed by any compact minimal
surfaces (horizons). Then
\begin{equation}\label{1.5}
m_{B}(\Omega) = \inf \{m_{ADM}(\widetilde M): (\widetilde M, g) \ {\rm is \ an \ admissible \
extension \ of} \ (\Omega, g)\},
\end{equation}
where $m_{ADM}(\widetilde M)$ is the ADM mass of $(\widetilde M, g)$, cf.~\cite{B1}. In \cite{HI} Huisken and Ilmanen 
have proved a number of basic properties of $m_{B}(\Omega)$, in particular that $m_{B}(\Omega)>0$ unless $(\Omega,g)$ is locally isometric to Euclidean space.
In \cite{Br} Bray discusses a similar definition, where the boundary $\partial \Omega$ is 
required to be outer-minimizing in $(\widetilde M, g)$. As will be seen below, the outer-minimizing 
property plays an important role in this paper, although for somewhat different reasons than in \cite{Br}. 

  Conjecturally, an extension $(\widetilde M, g)$ realizing the infimum in \eqref{1.5} is a solution
to the static vacuum Einstein equations \eqref{1.1} on $M = \widetilde M \setminus \Omega$ which
is Lipschitz, (but not smooth), across the junction $\partial \Omega$ and for which the
induced metric and mean curvature at the boundary of the interior and exterior regions agree:
$$g|_{\partial M} = g|_{\partial \Omega}, \ \ H_{\partial M} = H_{\partial \Omega},$$
leading to the boundary data \eqref{1.4}. Observe that the boundary data $(\gamma, H)$ have
the character of a mixed Dirichlet-Neumann type boundary value problem for the static
equations \eqref{1.1}, but the potential function $u$ is absent from the boundary data.
We point out that more standard Dirichlet or Neumann boundary data are not suitable for the
(static) Einstein equations, cf.~\cite{A3}.

\medskip

  In this paper, we develop a general framework for the Bartnik conjecture and make partial progress
on its resolution. To describe the setting, let ${\mathcal E}_{S} = {\mathcal E}_{S}^{m,\alpha}$
be the moduli space of AF static vacuum solutions $(M, g, u)$ on a given 3-manifold $M$ which
are $C^{m,\alpha}$ up to $\partial M$, $m \geq 3$.  The exact definition is given in \S 2, but
basically ${\mathcal E}_{S}$ is the space of all AF static vacuum metrics on $M$ modulo the
action of the group ${\rm Diff}_{1}$ of diffeomorphisms on $M$ equal to the identity on
$\partial M$.  Next, let $Met^{m,\alpha}(\partial M)$ be the space of
$C^{m,\alpha}$ metrics on $\partial M \simeq S^{2}$ and $C^{m-1,\alpha}(\partial M)$
be the space of $C^{m-1,\alpha}$ functions on $\partial M$. One thus has a natural
map, mapping a static vacuum solution to its Bartnik boundary data:
\begin{equation}\label{1.7}
\Pi_{B}:  {\mathcal E}_{S}^{m,\alpha} \rightarrow Met^{m,\alpha}(\partial M)\times
C^{m-1,\alpha}(\partial M),
\end{equation}
$$\Pi_{B}(g) = (\gamma, H).$$

\begin{theorem}\label{t1.1}
The space ${\mathcal E}_{S}^{m,\alpha}$ is a smooth (infinite dimensional) Banach manifold, and
the map $\Pi_{B}$ is $C^{\infty}$ smooth and Fredholm, of Fredholm index 0.
\end{theorem}

   Theorem 1.1 essentially amounts to the statement that the static vacuum Einstein equations
\eqref{1.1} with boundary conditions \eqref{1.4} form an elliptic boundary value problem,
modulo gauge transformations, i.e.~diffeomorphisms, and that one has a well-behaved local
existence theory for this problem. We note that this boundary value problem also has a variational
characterization, cf.~Proposition 3.7.

  Let ${\mathcal E}^{+}$ be the open submanifold of ${\mathcal E}_{S}$ for which
the mean curvature $H$ is positive, i.e.
$${\mathcal E}^{+} = (\Pi_{B})^{-1}(Met^{m,\alpha}(\partial M)\times
C_{+}^{m-1,\alpha}(\partial M)).$$
The Bartnik conjecture above may thus be rephrased to state that the smooth map $\Pi_{B}$,
restricted to the open submanifold ${\mathcal E}^{+}$,
\begin{equation}\label{1.8}
\Pi_{B}:  {\mathcal E}^{+} \rightarrow Met^{m,\alpha}(\partial M)\times
C_{+}^{m-1,\alpha}(\partial M),
\end{equation}
is surjective and injective, and hence, via the inverse function theorem, is a smooth diffeomorphism.

   However, this most optimistic version of the conjecture does not hold, in that $\Pi_{B}$ in \eqref{1.8} 
cannot be a diffeomorphism. To illustrate the problem consider (for example) the flat solution 
$g_{flat}$ with $u = 1$; there are boundaries $\partial M = S^{2}\subset {\mathbb R}^{3}$ 
given by embedded spheres $(S^{2}, \gamma_{i}, H_{i})$ with $H_{i}$ uniformly positive, which 
converge smoothly in ${\mathbb R}^{3}$ to a limit which is an immersed but not embedded sphere in 
${\mathbb R}^{3}$. Such a limit is then at the boundary $\partial {\mathcal E}^{+}$, but the limit 
boundary data $(\gamma, H) \in Met^{m,\alpha}(\partial M)\times C_{+}^{m-1,\alpha}(\partial M)$. 
In other words, the condition that the boundary data $(\gamma, H)$ is uniformly controlled in the 
target space is not sufficient to ensure that one stays within the class of manifolds-with-boundary.
(In particular, there cannot be a smooth inverse map to $\Pi_{B}$). As a concrete example, let 
$T^{2}$ be a torus of revolution embedded in ${\mathbb R}^{3}$
with $H > 0$. One may remove a (small) essential annulus from $T^{2}$ and smoothly
attach two embedded discs to obtain a 2-sphere $S^{2}$ with $H > 0$, cf.~Figure 1. This surface
may be deformed to obtain a curve $(S^{2})_{t}$, $t \in [0,1]$, of positive mean
curvature spheres which for $t < \frac{1}{2}$ are embedded and for $t \geq \frac{1}{2}$
are immersed, with a single self-intersection point of the discs at $t = \frac{1}{2}$.
(The same situation holds with any background static vacuum solution and varying
boundary $\partial M$ within $M$). This passage from embedded to immersed
behavior also shows that the boundary map $\Pi_{B}$ on ${\mathcal E}^{+}$ is not
proper.

\begin{figure}[htb]
\centering
\includegraphics[width=\hsize]{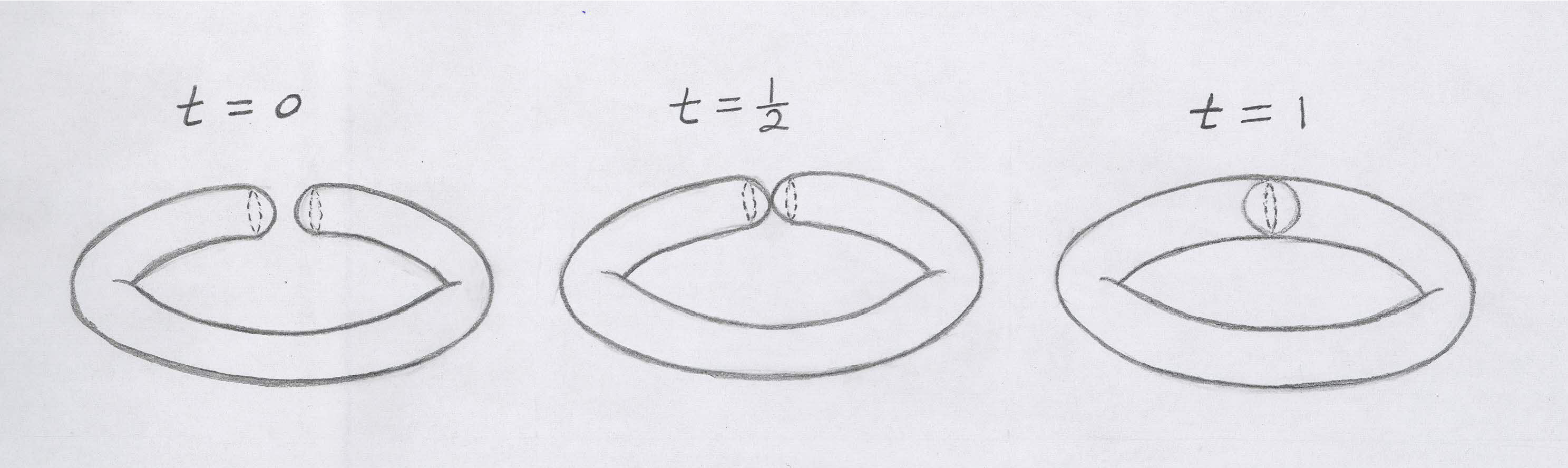}
\caption{\label{figure1} An illustration of the 1-parameter family of spheres $(S^{2})_{t}$, $t\in[0,1]$, in the process of passing from
embedding to immersion.}
\end{figure}

  The basic issue is in fact that of finding domains within ${\mathcal E}^{+}$
on which $\Pi_{B}$ is proper. Recall that a map between two topological spaces is proper,
if the preimage of any compact set is itself compact. In the current setting, $\Pi_{B}$ is proper on a domain ${\mathcal U}
\subset {\mathcal E}^{+}$ if whenever $(\gamma_{i}, H_{i})$ is a sequence of boundary
data converging to limit data $(\gamma, H)$ and $(M, g_{i}, u_{i})$ are any solutions
with $\Pi_{B}(g_{i}, u_{i}) = (\gamma_{i}, H_{i})$, then $(M, g_{i}, u_{i})$ converges, in a
subsequence, to a limit solution $(M, g, u) \in {\mathcal U}$. Here convergence is in the
topology of the target and domain spaces respectively. In other words, control
of the boundary data $(\gamma, H) = \Pi_{B}(g,u)$ implies global control of the solution
$(M, g, u)$ within ${\mathcal U}$. Equivalently, one needs apriori estimates
controlling the behavior of the full solution $(M, g, u)$ in terms of the boundary data
$(\gamma, H)$.

   Now let ${\mathcal E}^{o} \subset {\mathcal E}^{+}$ be the domain for which the
boundary $\partial M$ is strictly {\it outer-minimizing}, i.e.~for which
\begin{equation}\label{1.10}
area (\Sigma) > area (\partial M),
\end{equation}
for any surface $\Sigma \subset M$ homologous to $\partial M$ with $\Sigma \neq \partial M$.
(We point out that the examples in Figure 1 are {\em not} strictly outer-minimizing for $t$ sufficiently close to 
$\frac{1}{2}$). Clearly ${\mathcal E}^{o}$ is an open submanifold of ${\mathcal E}^{+}$. It is not known
(although likely to be true) that ${\mathcal E}^{o}$ is connected. Throughout the following, we
thus assume that ${\mathcal E}^{o}$ is taken to be the connected component containing the
standard flat exterior solution where $M$ is the exterior of the standard unit ball in
${\mathbb R}^{3}$, with $\partial M$ the round $S^{2}$ of radius 1. One then has a
natural boundary map
\begin{equation}\label{1.11}
\Pi^{o}: {\mathcal E}^{o} \rightarrow Met^{m,\alpha}(\partial M)\times
C_{+}^{m-1,\alpha}(\partial M).
\end{equation}

   The second main result of this paper is the following:
\begin{theorem}\label{t1.2}
The boundary map $\Pi^{o}$ in \eqref{1.11} is "almost" proper, in the following
sense. If $(\gamma_{i}, H_{i})$ is a sequence of boundary data converging to limit
data $(\gamma, H)$ in $Met^{m,\alpha}(\partial M)\times
C_{+}^{m-1,\alpha}(\partial M)$, and $(M, g_{i}, u_{i})\in {\mathcal E}^{o}$  are any
solutions with $\Pi^{o}(g_{i}, u_{i}) = (\gamma_{i}, H_{i})$, then $(M, g_{i}, u_{i})$ converges
in ${\mathcal E}_{S}^{m,\alpha}$, in a subsequence, to a limit solution $(M, g, u) \in
{\mathcal E}_{S}^{m,\alpha}$ which is (at least) weakly outer-minimizing, i.e.
\begin{equation}\label{1.12}
area (\Sigma) \geq area (\partial M),
\end{equation}
for $\Sigma$ as in \eqref{1.10}.
\end{theorem}

  Roughly speaking, Theorem 1.3 thus shows that static vacuum solutions $(M, g, u)$ with outer-minimizing 
boundary are controlled by their boundary data $(\gamma, H)$. The issue remains however of how to 
determine from the boundary data $(\gamma, H)$ whether the boundary is outer-minimizing in $(M, g, u)$. 
We point out that the full global property \eqref{1.10} is not actually necessary; Theorem
1.2 remains valid if \eqref{1.10} holds only in an arbitrarily small neighborhood of $\partial M$,
(depending on $(M, g, u)$), cf.~Remark 4.4.

\medskip 

    A smooth proper Fredholm map $F: B_{1} \rightarrow B_{2}$ between connected Banach manifolds
has a ${\mathbb Z}_{2}$-valued degree $deg_{{\mathbb Z}_{2}}F$, the Smale degree, cf.~\cite{Sm}. When the index of $F$
is zero, the degree is given by the number of preimages of a regular value modulo 2.
If the spaces or map have a suitable orientation, this can be extended to a ${\mathbb Z}$-valued
degree $deg_{\mathbb Z}F$, cf.~\cite{ET} for instance. Essentially an immediate consequence
of its definition and the Sard-Smale theorem \cite{Sm} is that if such a degree is non-zero, 
$$deg_{{\mathbb Z}_{2}}F \neq 0,$$
then $F$ is surjective. (Since the preimage of {\it any} regular value is non-emtpy). 

   The definition of degree above may be extended to maps which are almost proper
in the sense above, cf.~\cite{BFP} for instance. Thus, let $\partial {\mathcal E}^{o}$ be the
boundary of ${\mathcal E}^{o}$ within the space ${\mathcal E}_{S}$ of static vacuum solutions.
This is the space of solutions in ${\mathcal E}_{S}$ satisfying \eqref{1.12} but not \eqref{1.10}.
Let $Z = \Pi_{B}(\partial {\mathcal E}^{o}) \subset Met^{m,\alpha}(\partial M)\times
C^{m-1,\alpha}(\partial M)$ be the image of $\partial {\mathcal E}^{o}$ under the boundary
map $\Pi_{B}$ and let ${\mathcal E}^{P} = (\Pi^{o})^{-1}([Met^{m,\alpha}(\partial M)\times
C^{m-1,\alpha}(\partial M)]\setminus Z)$ be the corresponding inverse image. Then, as
discussed in \S 5, the induced boundary map (restriction of $\Pi^{o}$ to ${\mathcal E}^{P}$)
$$\Pi^{P}: {\mathcal E}^{P} \rightarrow [Met^{m,\alpha}(\partial M)\times
C^{m-1,\alpha}(\partial M)]\setminus Z,$$
is proper. In particular, ${\mathcal E}^{P}$ has a finite number of connected components
${\mathcal E}^{P_{i}}$ and the induced boundary map $\Pi^{P_{i}}$ on ${\mathcal E}^{P_{i}}$
has a well-defined ${\mathbb Z}_{2}$-valued degree (with respect to a component of the
target space). We also expect that $\Pi^{P_{i}}$ has a well-defined
${\mathbb Z}$-valued degree.

   Let ${\mathcal E}^{P_{0}}$ be the component of ${\mathcal E}^{P}$ containing the standard round
exterior flat solution as following \eqref{1.10}, and let ${\mathcal T}_{0}$ be the component of
$[Met^{m,\alpha}(\partial M)\times C_{+}^{m-1,\alpha}(\partial M)] \setminus Z$ containing
the corresponding standard boundary data $(\gamma_{+1}, 2)$. One thus has the boundary map
\begin{equation}\label{1.13}
\Pi^{P_{0}}: {\mathcal E}^{P_{0}} \rightarrow {\mathcal T}_{0}.
\end{equation}
A further main result of this paper is the computation of this degree:
\begin{theorem}\label{t1.3}
For the map $\Pi^{P_{0}}$ in \eqref{1.13}, one has
\begin{equation}\label{1.14}
deg_{{\mathbb Z}_{2}}\Pi^{P_{0}} = 1.
\end{equation}
\end{theorem}
The proof of Theorem 1.3 is based on the well-known black hole uniqueness theorem for the
Schwarzschild metric, cf.~\cite{I}, \cite{R}, \cite{BM}.

   It follows that the boundary map $\Pi^{P_{0}}$ maps ${\mathcal E}^{P_{0}}$ surjectively
onto the component ${\mathcal T}_{0}$ of the target $[Met^{m,\alpha}(\partial M)\times
C_{+}^{m-1,\alpha}(\partial M)] \setminus Z$. In particular, the image of $\Pi^{P_{0}}$
and so the image of $\Pi_{B}$ has non-empty interior in the target space
$Met^{m,\alpha}(\partial M)\times C_{+}^{m-1,\alpha}(\partial M)$; this has
not been previously known, cf.~Remark 5.4 for further discussion.

\medskip 

   Of course, the main issue at this point is what can be said about the structure of the
set $Z$? It is of first category (so non-generic) but its more detailed structure awaits future study.
An alternate approach bypassing the issue of the boundary $\partial {\mathcal E}^{o}$
would be to find conditions on the boundary data $(\gamma, H)$ which imply $\partial M$
is outer-minimizing as in \eqref{1.10} in any static vacuum extension $(M, g, u)$ of
$(\gamma, H)$. For example, in ${\mathbb R}^{3}$ (where $u = 1$) convexity suffices,
which is expressed in terms of boundary data as $K_{\gamma} > 0$ where $K$ is the
Gauss curvature. It is an open problem whether this condition also suffices for general
static vacuum solutions.

\smallskip

  The contents of the paper are briefly as follows. In \S 2, we present background
information on the structure of static vacuum solutions and choices of gauge.
Section 3 discusses elliptic boundary value problems for the Einstein equations
and proves the basic structure theorems on the moduli space ${\mathcal E}_{S}$
and the boundary map $\Pi_{B}$, including Theorem 1.1.  In \S 4, we then prove the requisite
a priori estimates and establish the almost properness of $\Pi^{o}$ on ${\mathcal E}^{o}$,
proving Theorem 1.2. Finally, \S 5 contains the computation of the degree of $\Pi^{P_{0}}$
and closes with several related remarks.

\smallskip

   We thank Robert Bartnik, Piotr Chru\'sciel, Gerhard Huisken and Xin Zhou for their
interest and comments on this work. We are especially grateful to Simon Brendle for
pointing out an error in a previous version of the paper.

\section{Background Discussion}

\setcounter{equation}{0}

  Let $M$ be a 3-manifold with compact boundary $\partial M$, and with a single
open end $E$. (All of the results of this section and of \S 3
hold in all dimensions, but for simplicity, we restrict to dimension 3).
A priori, $\partial M$ need not be connected, although this will be assumed
later on. As following \eqref{1.2}-\eqref{1.3}, we let $\mathcal{M} = S^{1}\times M$. Almost all of 
the discussion and computation in Sections 2 and 3 is carried out on the $4$-manifold $\mathcal{M}$ 
and $g_{\mathcal M}$ will often be denoted $g$ for notational simplicity. 

  Let $Met_{S}(\mathcal{M}) = Met_{S}^{m,\alpha}(\mathcal{M})$ be the space of complete 
(up to the boundary) $C^{m,\alpha}$ static metrics on $\mathcal{M}$, i.e.~metrics of the form 
\eqref{1.2}, $m \geq 2$. One has
\begin{equation}\label{2.1}
Met_{S}^{m,\alpha}(\mathcal{M}) \simeq Met^{m,\alpha}(M)\times C_{+}^{m,\alpha}(M),
\end{equation}
where $C_{+}^{m,\alpha}(M)$ is the space of positive $C^{m,\alpha}$ functions on
$M$. The space ${\mathbb E}_{S} = {\mathbb E}_{S}(\mathcal{M})$ of static Einstein (Ricci-flat)
metrics on $\mathcal{M}$ is equivalent to the space of pairs $g_{\mathcal{M}} = (g_{M}, u) \in
Met(M)\times C_{+}(M)$ satisfying \eqref{1.1} or \eqref{1.3} (the smoothness
indices will be occasionally dropped when unimportant). It is well-known \cite{M} that away 
from the boundary, solutions of the static vacuum equations are analytic in appropriate coordinates.

  Recall that a complete metric $g \in Met^{m,\alpha}(E)$ on an end $E$ is
asymptotically flat if $E$ is diffeomorphic to ${\mathbb R}^{3}\setminus B$,
where $B$ is a 3-ball, and there exists a diffeomorphism $F:
{\mathbb R}^{3}\setminus B \rightarrow E$ such that, in the chart $F$,
\begin{equation}\label{2.2}
g_{ij} = \delta_{ij} + O(r^{-1}), \ \ \partial_{k}g_{ij} = O(r^{-2}), \ \
\partial_{k}\partial_{\ell}g_{ij} = O(r^{-3}),
\end{equation}
in the standard Euclidean coordinates on ${\mathbb R}^{3}$. The static vacuum
equations \eqref{1.1} are invariant under multiplication of the potential $u$ by
constants. Throughout the paper, we assume that $u$ is normalized so that
$u \rightarrow 1$ at infinity, and that $u$ is asymptotically constant in
the sense that
\begin{equation}\label{2.3}
u = 1 + O(r^{-1}), \ \ \partial_{k}u = O(r^{-2}), \ \
\partial_{k}\partial_{\ell}u = O(r^{-3}).
\end{equation}
Thus the 4-metric $g_{\mathcal{M}}$ is asymptotic to the product $S^{1}\times {\mathbb R}^{3}$
at infinity.

  It is proved in \cite{A2} that ends of static vacuum solutions $(M, g, u)$ are
either asymptotically flat or parabolic, where parabolic is understood in the sense
of potential theory; equivalently, parabolic ends have small volume growth
in that the area of geodesic spheres grows slower than $r^{1+\varepsilon}$,
for any $\varepsilon > 0$. Moreover, asymptotically flat ends are strongly
asymptotically flat in that the metric and potential have asymptotic
expansions of the form
\begin{equation}\label{2.4}
g_{ij} = \left(1 + \frac{2m}{r}\right)\delta_{ij} + \cdots, \ \
u = 1 - \frac{2m}{r} + \cdots,
\end{equation}
where the mass $m$ may apriori be any value $m \in {\mathbb R}$, cf.~also \cite{KO}. These two behaviors,
asymptotically flat and parabolic, are radically different and there is no curve
of asymptotic structures for static vacuum solutions which joins them. The finer
behavior of asymptotically flat ends is described by the mass parameter $m$
and higher multipole moments, cf.~\cite{BS}.

\medskip

 Let $\widetilde g$ be a fixed asymptotically flat (background) metric in
${\mathbb E}_{S}$; henceforth ${\mathbb E}_{S}$ will denote the space of
asymptotically flat static vacuum Einstein solutions. The static Einstein
equations are not elliptic, due to their invariance under diffeomorphisms,
and for several reasons one needs to choose an elliptic gauge. To begin,
we consider the Bianchi gauge, and define
\begin{equation} \label{2.5}
\Phi_{\widetilde g}: Met_{S}^{m,\alpha}(\mathcal{M}) \rightarrow  S^{m-2,\alpha}(\mathcal{M}),
\end{equation}
$$\Phi_{\widetilde g}(g) = Ric_{g} + \delta_{g}^{*}\beta_{\widetilde g}(g),$$
where $\beta_{\widetilde g}$ is the Bianchi operator, $\beta_{\widetilde g}(g)
= \delta_{\widetilde g}(g) + \frac{1}{2}dtr_{\widetilde g}(g)$. Also,
$(\delta^{*}X)(A,B) = \frac{1}{2}(\langle \nabla_{A}X, B\rangle  + \langle
\nabla_{B}X, A\rangle )$ and $\delta X = -tr(\delta^{*}X)$ is the divergence.
The operator $\Phi_{\widetilde g}$ is a $C^{\infty}$ smooth map into the
space $S^{m-2,\alpha}(\mathcal{M})$ of static symmetric bilinear forms on $\mathcal{M}$; note that here a symmetric bilinear form
is referred to as static if its components do not depend on time and the mixed time/space components vanish.

   Using standard formulas for the linearization of the Ricci and scalar
curvatures, cf.~\cite{Be} page 63 for instance, the linearization of $\Phi$ at $g =
\widetilde g \in {\mathbb E}_{S}$ is given by
\begin{equation} \label{2.6}
L(h) = 2(D\Phi_{g})(h) = D^{*}Dh - 2R(h),
\end{equation}
where $R(h)(X,Y) = \langle R(e_{i}, X)Y, h(e_{i})\rangle$ with $e_{i}$ an orthonormal basis. 
Clearly, $L$ is elliptic and formally self-adjoint. In \S 3 we will discuss boundary value problems 
for $\Phi$ and $L$.

  Next, the asymptotic behavior in the asymptotically flat end $E$ requires
the introduction of suitable weighted function spaces. We will use the
standard weighted H\"older spaces, although one could equally well use
weighted Sobolev spaces. Thus, define $Met_{\delta}^{m,\alpha}(\mathcal{M}) \subset
Met_{S}^{m,\alpha}(\mathcal{M})$ to be the subspace of metrics which decay to Euclidean
data at a rate $r^{-\delta}$ at infinity; more precisely, the component
functions $g_{ij}$ and $u$ of $g_{\mathcal{M}}$ should satisfy
$$g_{ij} - \delta_{ij} \in C_{\delta}^{m,\alpha}({\mathbb R}^{3}\setminus B),
\ \ u - 1 \in C_{\delta}^{m,\alpha}({\mathbb R}^{3}\setminus B).$$
Here $C_{\delta}^{m}$ consists of functions $v$ such that
$$||v||_{C_{\delta}^{m}} = \sum_{k=0}^{m}\sup r^{k+\delta}|\nabla^{k}v| <
\infty,$$
while $C_{\delta}^{m,\alpha}$ consists of functions such that
$$||v||_{C_{\delta}^{m,\alpha}} = ||v||_{C_{\delta}^{m}} + \sup_{x,y}
[\min(r(x),r(y))^{m+\alpha+\delta}\frac{|\nabla^{m}v(x) - \nabla^{m}v(y)|}
{|x - y|^{\alpha}}] < \infty,$$
cf.~\cite{B1}, \cite{LP}. Throughout the following, we assume the decay rate
$\delta$ is fixed, and chosen to satisfy
\begin{equation}\label{2.7}
{\tfrac{1}{2}} < \delta < 1.
\end{equation}
It is well-known, cf.~\cite{B1}, \cite{LP}, that the Laplacian on functions, and
Laplace-type operators on tensors, as in \eqref{2.6}, are Fredholm when acting
on these weighted H\"older spaces.

  The map $\Phi$ in \eqref{2.5} clearly induces a smooth map
\begin{equation} \label{2.8}
\Phi : Met_{\delta}^{m,\alpha}(\mathcal{M}) \rightarrow S_{\delta}^{m-2,\alpha}(\mathcal{M}).
\end{equation}
Observe that $g$ is Einstein if $\Phi_{\widetilde g}(g) = 0$ and $\beta_{\widetilde g}(g)
= 0$, so that $g$ is in Bianchi gauge with respect to $\widetilde g$. (Although
$\Phi_{\widetilde g}$ is defined for all $g\in Met_{\delta}^{m,\alpha}(\mathcal{M})$, we will
only consider it acting on $g$ near $\widetilde g$).

 Given $\widetilde g \in {\mathbb E}_{S}$, let $Met_{C}^{m,\alpha}(\mathcal{M}) \subset
Met_{\delta}^{m,\alpha}(\mathcal{M})$ be the space of $C^{m,\alpha}$ smooth static AF
Riemannian metrics on $\mathcal{M}$, satisfying the Bianchi gauge constraint
\begin{equation}\label{2.9}
\beta_{\widetilde g}(g) = 0 \ \ {\rm at} \ \ \partial \mathcal{M}.
\end{equation}
As above,
$$\Phi : Met_{C}^{m,\alpha}(\mathcal{M}) \rightarrow S_{\delta}^{m-2,\alpha}(\mathcal{M}).$$
Similarly let $Z_{C}^{m,\alpha}$ be the space of metrics $g \in
Met_{C}^{m,\alpha}(\mathcal{M})$ satisfying $\Phi_{\widetilde g}(g) = 0$, and let
\begin{equation} \label{2.10}
{\mathbb E}_{C} \subset  Z_{C}
\end{equation}
be the subset of static Einstein metrics $g = g_{\mathcal{M}}$, $Ric_{g} = 0$ in $Z_{C}$.
The following result justifies the use of the operator $\Phi$ to study
${\mathbb E}_{C}$.

\begin{proposition}\label{p2.1}
Any static metric $g = g_{\mathcal{M}} \in Z_{C}$ sufficiently close to $\widetilde g$
is necessarily Einstein, $g \in {\mathbb E}_{C}$. Moreover, this also holds
infinitesimally in the following sense. Let $\kappa$ be an infinitesimal
deformation of $g\in Z_{C}$, i.e.~$\kappa\in Ker D\Phi$.
If $\beta_{\widetilde g}(g) = 0$, (for example $\widetilde g = g$),
then
\begin{equation} \label{2.11}
\beta_{\widetilde g}(\kappa) = 0,
\end{equation}
and $\kappa$ is an infinitesimal Einstein deformation, i.e.~the
variation of $g$ in the direction $\kappa$ preserves \eqref{1.3} to
$1^{\rm st}$ order.
\end{proposition}

{\bf  Proof:} Since $g\in Z_{C}$, one has $\Phi (g) = 0$, i.e.
$$Ric_{g} + \delta_{g}^{*}\beta_{\widetilde g}(g) = 0.$$
Applying the Bianchi operator $\beta_{g}$ and using the Bianchi
identity $\beta_{g}(Ric_{g}) = 0$ gives
\begin{equation} \label{2.12}
\beta_{g}(\delta_{g}^{*}(\beta_{\widetilde g}(g))) = 0.
\end{equation}
Set $V = \beta_{\widetilde g}(g)$, and notice that a simple computation produces the Weitzenbock formula
$2\beta_{g}\delta_{g}^{*}(V) = D^{*}DV - Ric(V)$. Also,
since $g, \widetilde g \in Met_{\delta}^{m,\alpha}(\mathcal{M})$, $V \in
\chi_{1+\delta}^{m-1,\alpha}(M)$, where $\chi_{1+\delta}^{m-1,\alpha}(M)$
is the space of vector fields whose components are in
$C_{1+\delta}^{m-1,\alpha}(M)$. When acting on vector fields $V$ with
$V = 0$ on $\partial M$, as in \eqref{2.9}, the operator $D^{*}D$ is positive,
with trivial kernel. Namely, if $W \in C_{1+\delta}^{m-1,\alpha}$ is in
the kernel of $D^{*}D$, then integrating by parts gives
$$\int_{B(r)}|DW|^{2} + \int_{S(r)}\langle W, \nabla_{N}W \rangle = 0,$$
where $B(r) = \{x \in M: dist(x, \partial M) \leq r\}$ and $N$ is the outward
unit normal. (Since $W = 0$ on $\partial M$, there is no boundary term at
$\partial M$). Letting $r \rightarrow \infty$, the boundary integral tends to
0 and so $DW = 0$, which in turn implies $W = 0$.

  Since $D^{*}D$ is self-adjoint and Fredholm, it has a smallest positive
eigenvalue bounded away from 0. For $g$ sufficiently close to $\widetilde g$,
$|Ric|\sim 0$ pointwise and $Ric(V) \in C_{3+\delta}^{m-3,\alpha}(M)$,
so we may assume that $2\beta_{g}\delta_{g}^{*}$ is a positive operator
on $V$. Hence, again since $V = 0$ on $\partial M$, the only solution of
\eqref{2.12} is $V = 0$, which implies $g \in {\mathbb E}_{C}$.

  To prove the second statement, let $g_{t} = g+t\kappa$. Applying the
Bianchi operator $\beta_{g_{t}}$ to $\Phi(g_{t})$ gives
\begin{equation}\label{2.13}
\beta_{g_{t}}\Phi (g_{t}) =
\beta_{g_{t}}\delta_{g_{t}}^{*}(\beta_{\widetilde g}(g_{t})).
\end{equation}
Taking the derivative with respect to $t$ at $t = 0$, one has
$(\beta_{g_{t}}\Phi (g_{t}))'  = \beta'\Phi  + \beta\Phi'$. Both
terms here vanish since $g\in Z_{C}$ and $\kappa$ is formally tangent
to $Z_{C}$. Hence the variation of the right hand side of \eqref{2.13}
vanishes. Since $\beta_{\widetilde g}(g) = 0$, this gives
$\beta_{g}\delta_{g}^{*}(\beta_{\widetilde g}(\kappa)) = 0$. The
equation \eqref{2.11} then follows exactly as following \eqref{2.12},
with $V = \beta_{\widetilde g}(\kappa)$.

{\endproof}

 Let ${\mathcal D}_{1}^{m+1,\alpha}$ denote the space of $C_{\delta}^{m+1,\alpha}$
static diffeomorphisms of $\mathcal{M}$ which equal the identity on $\partial \mathcal{M}$. These are
diffeomorphisms which decay to the identity at the rate $r^{-\delta}$ and are independent
of the $t$ or $\theta$-variable in \eqref{1.2}. The group ${\mathcal D}_{1}$ acts
freely and continuously on $Met(\mathcal{M})$ and ${\mathbb E}_{S}$ by pullback and one has
the following local slice theorem for this action; we refer to \cite{A3} for
the proof.

\begin{lemma}\label{l2.2}
Given any $\widetilde g \in {\mathbb E}_{S}^{m,\alpha}$ and $g \in Met_{\delta}^{m,\alpha}(\mathcal{M})$
near $\widetilde g$, there exists a unique diffeomorphism $\phi \in
{\mathcal D}_{1}^{m+1,\alpha}$ close to the identity, such that
\begin{equation} \label{2.14}
\beta_{\widetilde g}(\phi^{*}g) = 0.
\end{equation}
In particular, $\phi^{*}g\in Met_{C}^{m,\alpha}(\mathcal{M})$.
\end{lemma}

  Lemma 2.2 implies that if $g \in {\mathbb E}_{S}^{m,\alpha}$ is a static Einstein
metric near $\widetilde g$, then $g$ is isometric, by a diffeomorphism in
${\mathcal D}_{1}^{m+1,\alpha}$, to an Einstein metric in ${\mathbb E}_{C}^{m,\alpha}$,
so that ${\mathbb E}_{C}^{m,\alpha}$ is a slice for ${\mathbb E}_{S}^{m,\alpha}$ under
the action of ${\mathcal D}_{1}^{m+1,\alpha}$.

   To prove that the moduli space ${\mathcal E}$ is a smooth Banach manifold,
(cf.~Theorem 3.6), it is important to have a gauge with choice of boundary
data in which the Einstein equations form a self-adjoint elliptic boundary value problem.
This is not the case for the operator $\Phi$ and we are not aware of geometrically natural
self-adjoint boundary conditions for $\Phi$. For this reason, we will also consider
another natural gauge, namely the divergence-free gauge.

   To do this, instead of $\Phi$, consider
\begin{equation}\label{2.15}
\hat \Phi(g) = \hat \Phi_{\widetilde g}(g) = Ric_{g} - \frac{s}{2}g +
\delta_{g}^{*}\delta_{\widetilde g}g,
\end{equation}
where $s$ is the scalar curvature of $g = g_{\mathcal{M}}$. The linearization of $\hat \Phi$
at $g = \widetilde g \in {\mathbb E}_{S}$ is given by
\begin{equation} \label{2.16}
\hat L(h) = 2(D\hat \Phi_{\widetilde g})_{g}(h) = D^{*}Dh - 2R(h) - (D^{2}trh +
\delta\delta h \,g) + \Delta trh \,g .
\end{equation}

  In analogy to \eqref{2.9}, define
\begin{equation}\label{e2.17}
Met_{D}^{m,\alpha}(M) = \{g \in Met_{\delta}^{m,\alpha}(\mathcal{M}): \delta_{\widetilde g}g = 0
\ \ {\rm at} \ \ \partial \mathcal{M}\}.
\end{equation}
Similarly, let $Z_{D}^{m,\alpha} = \hat \Phi^{-1}(0)\cap Met_{D}^{m,\alpha}(M)$ and
${\mathbb E}_{D} \subset Z_{D}$ be the space of Einstein metrics in divergence-free
gauge with respect to $\widetilde g \in {\mathbb E}_{S}$.

   It is easy to see that Proposition 2.1 and Lemma 2.2 hold in this divergence-free
gauge in place of the previous Bianchi gauge, with essentially the same proof. Thus
${\mathbb E}_{D} = Z_{D}$ and for $g \in {\mathbb E}_{D}$,
\begin{equation}\label{2.18}
\delta_{\widetilde g}g = 0 \ \ {\rm on} \ \ \mathcal{M}.
\end{equation}
Moreover, the diffeomorphism group ${\mathcal D}_{1}$ transforms one gauge choice
uniquely to the other. For instance, suppose $\beta_{\widetilde g}(g) = 0$. Then
we claim there is a unique $\phi \in {\mathcal D}_{1}^{m+1,\alpha}$ such that
\begin{equation}\label{2.19}
\delta_{\widetilde g}(\phi^{*}g) = 0.
\end{equation}
At the linearized level, with $g = \widetilde g$, this amounts to finding a vector
field $V$ with $V = 0$ on $\partial \mathcal{M}$ such that if $\beta_{\widetilde g}h = 0$ then
$\delta_{\widetilde g}(h + \delta^{*}V) = 0$. This equation is equivalent to the
equation $\delta \delta^{*}V = \frac{1}{2}dtr h$, which is uniquely solvable for
$V$ with $V = 0$ on $\partial \mathcal{M}$. The local result in \eqref{2.19} then follows
from the inverse function theorem.

\section{The Moduli Space}

\setcounter{equation}{0}

  In this section, we study boundary value problems for the elliptic operators
$\Phi$ and $\hat \Phi$, and use this to prove that the moduli space ${\mathcal E}_{S}$
of static vacuum solutions is a smooth Banach manifold for which the boundary map
$\Pi_{B}$ is Fredholm, of Fredholm index 0, cf.~Theorem 3.6.

   We begin with the Bianchi-gauged Einstein operator $\Phi$ in \eqref{2.5}, i.e.
$$\Phi_{\widetilde g}(g) = Ric_{g} + \delta_{g}^{*}\beta_{\widetilde g}(g).$$
Let $A$ denote the $2^{\rm nd}$ fundamental form of $\partial M$ in $M$,
$A(X,Y) = \langle \nabla_{X}N, Y \rangle$, where $N$ is the unit inward normal
into $M$, $X, Y$ tangent to $\partial M$. Similarly, let $H_{M} = tr A$ denote the
mean curvature of $\partial M$ in $M$. Throughout the paper $W^{T}$ will denote
the restriction or the orthogonal projection of a tensor $W$ to $T(\partial N)$
or $T(\partial M)$.

\begin{proposition}\label{p3.1}
Near any given background solution $\widetilde g \in {\mathbb E}_{S}^{m,\alpha}$, the
operator $\Phi = \Phi_{\widetilde g}$ in \eqref{2.5} with boundary conditions:
\begin{equation}\label{3.1}
\beta_{\widetilde g}(g) = 0, \ \ g|_{T(\partial M)} = \gamma_{M}, \ \
H_{M} = h_{M} \ \ \textit{at} \ \ \partial \mathcal{M},
\end{equation}
is an elliptic boundary value problem of Fredholm index 0.
\end{proposition}

  Here the induced metric $\gamma_{M}$ is in $Met^{m,\alpha}(\partial M)$ while the
mean curvature $h_{M}$ of $\partial M$ in $(M, g_{M})$ is in $C^{m-1,\alpha}(\partial M)$.
Note that the potential $u$ does not enter this boundary data and so is formally
undetermined at $\partial M$. Also the static property implies that
$\beta_{\widetilde g}(g)$ vanishes in the vertical direction,
$\beta_{\widetilde g}(g)(\partial_{\theta}) = 0$.

{\bf Proof:} It suffices to prove that the leading order part of the linearized
operators forms an elliptic system. Recall from \eqref{2.6} that the linearization of
$\Phi$ at $\widetilde g = g$ is given by
$$L(h) = 2(D\Phi_{g})(h) = D^{*}Dh - 2R(h).$$
The leading order symbol of $L = 2D\Phi$ at $\xi'$ is
\begin{equation}\label{3.2}
\sigma(L) = -|\xi'|^{2}I,
\end{equation}
where $I$ is the $Q\times Q$ identity matrix, with $Q = (n(n+1)/2) + 1$;
$Q$ is the sum of the dimension of the space of symmetric bilinear forms on ${\mathbb R}^{n}$,
together with the extra vertical $S^{1}$ direction. Here $n = 3$ but we give the
proof for general dimensions. For static metrics, all components of the metric are
locally functions on ${\mathbb R}^{n}$, and all derivatives in the vertical $S^{1}$
direction are trivial. In the following, the subscript 0 represents the direction
normal to $\partial M$ in $M$, (or $\partial \mathcal{M}$ in $\mathcal{M}$), subscript 1 denotes the
vertical direction, tangent to $S^{1}$, while indices $2$ through $n$ represent the
directions tangent to $\partial M$. Note that one has $h_{1\alpha} = 0$, for all $\alpha
\neq 1$. The positive roots of \eqref{3.2} are $i|\xi|$, where $\xi' = (\xi_{0}, \xi)$,
with multiplicity $Q$ at $\xi \in T^{*}({\mathbb R}^{n})$.

  Writing $\xi' = (z, \xi_{i})$, $i = 2, \dots , n$, (as above $\xi_{1} = 0$), the
symbols of the leading order terms in the boundary operators are given by:
$$-2izh_{0k} - 2i\sum_{j\geq 2} \xi_{j}h_{jk} + i\xi_{k}tr h = 0, \ \ k \geq 2,$$
$$-2izh_{00} - 2i\sum_{k\geq 2} \xi_{k}h_{0k} + iztr h = 0,$$
$$h^{T} = (\gamma')^{T}, \ \ {\rm and} \ \ (H_{M})'_{h} = \omega.$$
This gives $n + \frac{n(n-1)}{2} + 1 = Q$ boundary equations, as required.
Ellipticity requires that the operator defined by the boundary symbols above
has trivial kernel when $z$ is set to the root $i|\xi|$. Carrying this out
then gives the system
\begin{equation}\label{3.3}
2|\xi|h_{0k} - 2i\sum_{j\geq 2} \xi_{j}h_{jk} + i\xi_{k}tr h = 0, \ \ k \geq 2,
\end{equation}
\begin{equation}\label{3.4}
2|\xi|h_{00} - 2i\sum_{k\geq 2} \xi_{k}h_{0k} - |\xi|tr h = 0,
\end{equation}
\begin{equation}\label{3.5}
h_{11} = \phi, \ \ h^{T} = 0, \ \ (H_{M})'_{h} = 0,
\end{equation}
where $\phi$ is an undetermined function.

  Multiplying \eqref{3.3} by $i\xi_{k}$ and summing gives, via \eqref{3.5},
$$2|\xi|i\sum_{k\geq 2}\xi_{k}h_{0k} =  |\xi|^{2}tr h.$$
Substituting \eqref{3.4} on the term on the left above then gives
$$2|\xi|^{2}h_{00} - 2|\xi|^{2}tr h =  0.$$
Since $tr h = h_{00} + \phi$, it follows that $\phi = 0$.

  Next, to compute $H_{M}'$, we first observe that in general
\begin{equation}\label{3.6}
2A'_{h} = \nabla_{N}h + 2A\circ h - 2\delta^{*}(h(N)^{T}) - \delta^{*}(h_{00}N).
\end{equation}
This follows by differentiating the defining formula $2A = {\mathcal L}_{N}g$, and
using the identities $2N'_{h} = - 2h(N)^{T} - h_{00}N$, ${\mathcal L}_{N}h =
\nabla_{N}h + 2A\circ h$. Since $H_{M} = tr_{M}A$, $H'_{h} = tr_{M}A'_{h} - tr_{M}
A \circ h$ and so
\begin{equation}\label{3.7}
2(H_{M})'_{h} = tr_{M}(\nabla_{N}h - 2\delta^{*}(h(N)^{T}) - \delta^{*}(h_{00}N)).
\end{equation}
Hence the symbol of $2(H_{M})'_{h}$ is given by
$\sum_{k\geq 2} (iz h_{kk} - 2i\xi_{k}h_{0k})$. Setting this to
0 at the root $z = i|\xi|$ gives
\begin{equation}\label{3.8}
\sum_{k\geq 2}(|\xi| h_{kk} + 2i\xi_{k}h_{0k}) = 0.
\end{equation}
Via \eqref{3.5}, this gives $-2i\sum_{k\geq 2}\xi_{k}h_{0k} = 0$, and substituting
this in \eqref{3.4} and using the fact that $\phi = 0$ gives
$$2|\xi|h_{00} - |\xi|h_{00} = 0,$$
so that $h_{00} = 0$. It follows from \eqref{3.3} that $h_{0k} = 0$ and hence
$h = 0$. This proves ellipticity of the boundary value problem \eqref{3.1} and
the Fredholm property follows from the fact that the Laplace-type operator $L$
is Fredholm on $Met_{\delta}^{m,\alpha}$, cf.~\cite{LP}.

  Finally, it is straightforward to verify that the boundary data \eqref{3.1}
may be continuously deformed through elliptic boundary data to elliptic boundary
data for which $L$ is self-adjoint and so of index 0. This is proved in \cite{A3}
in a slightly different setting and the proof carries over here with only
minor change, and so we refer to \cite{A3} for further details. The homotopy
invariance of the index then completes the proof.

{\endproof}

  As noted in \S 2, we are not aware of a geometrically natural self-adjoint elliptic
boundary value problem for $\Phi$. In particular, the boundary conditions \eqref{3.1}
are not self-adjoint. This property is important for the proof of Theorem 3.6, and
for this reason, we turn to the operator $\hat \Phi$ in \eqref{2.15} with linearization
at $\widetilde g = g$ given by $\hat L$ in \eqref{2.16}.

  Regarding boundary conditions for $\hat L$, for $h \in S_{\delta}^{m-2,\alpha}(\mathcal{M})$,
let $h^{T} = h|_{\partial \mathcal{M}}$ and $[h^{T}]_{0}$ be the projection of $h^{T}$ onto
the space of forms trace-free with respect to $\gamma = \gamma_{\mathcal{M}}$. Similarly,
$H'_{h}$ denotes here the linearization of the mean curvature $H = H_{\mathcal{M}}$ of
$\partial \mathcal{M} \subset \mathcal{M}$.

  We then have:
\begin{lemma}\label{l3.2}
The operator $\hat L$ with boundary conditions
\begin{equation}\label{3.9}
\delta h = 0, \ \ [h^{T}]_{0} = 0, \ \ H'_{h} = 0,
\end{equation}
is a self-adjoint elliptic operator. Moreover, under the first two conditions
$\delta h = 0$ and $[h^{T}]_{0} = 0$, the operator $\hat L$ is self-adjoint exactly
for the boundary condition $H'_{h} = 0$.
\end{lemma}

{\bf Proof:} It is a rather long (and uninteresting) calculation to prove that the
operator $\hat L$ with boundary data \eqref{3.9} forms an elliptic system; this has been
verified by computer computation using Maple. More conceptually, instead we will make
use of Proposition 3.1 to simplify the proof. First, recall, \cite{ADN}, \cite{T},
that ellipticity of a boundary value problem is equivalent to the existence of a
uniform estimate
\begin{equation}\label{3.10}
||h||_{C^{m,\alpha}} \leq C(||\hat L(h)||_{C^{m-2,\alpha}} + ||B_{j}(h)||_{C^{m-j,\alpha}}
+ ||h||_{C^{0}}),
\end{equation}
where $B_{j}$ is the part of the boundary operator of order $j$, together with
such an estimate for the adjoint operator. As seen below, the boundary value
problem is self-adjoint, so it suffices to establish \eqref{3.10}.

  First, it is simple to prove \eqref{3.10} for $L$ in place of $\hat L$ via
a slight modification of the proof of Proposition 3.1. Namely, for the boundary
condition $[h^{T}]_{0} = 0$, we have $h^{T} = \phi \gamma$ on $\partial \mathcal{M}$, (in place
of \eqref{3.5}). Note also that \eqref{3.3}-\eqref{3.4} hold, but without the
$tr h$ terms. The analog of \eqref{3.3} then gives
$$|\xi|h_{0k} = i\xi_{k}\phi,$$
and hence, via the analog of \eqref{3.4},
$$|\xi|^{2}h_{00} = -|\xi|^{2}\phi,$$
so that $h_{00} + \phi = 0$. Next, via the condition $H'_{h} = 0$, the analog of
\eqref{3.8} becomes
$$\sum_{k \geq 1}(|\xi|h_{kk} + 2i\xi_{k}h_{0k}) = 0,$$
which gives
$$n|\xi|\phi = -2i\sum \xi_{k}h_{0k} = 2|\xi|\phi.$$
Since $n \geq 3$, this implies $\phi = 0$, and so $h_{00} = 0$, hence $h_{0k} = 0$.
It follows that $h = 0$, which proves ellipticity of $L$ with the boundary
conditions \eqref{3.9}. Thus, \eqref{3.10} holds with $L$ in place of $\hat L$.

  Next, one has
\begin{equation}\label{3.11a}
\hat L = L - (D^{2}tr h - \Delta tr h \,g) - \delta \delta h \,g.
\end{equation}
Thus to prove \eqref{3.10}, it suffices to prove
\begin{equation}\label{3.11}
||\delta h||_{C^{m-1,\alpha}} \leq C(||\hat L(h)||_{C^{m-2,\alpha}} + ||B_{j}(h)||_{C^{m-j,\alpha}}
+ ||h||_{C^{0}}),
\end{equation}
\begin{equation}\label{3.12}
||D^{2}tr h||_{C^{m-2,\alpha}} \leq C(||\hat L(h)||_{C^{m-2,\alpha}} + ||B_{j}(h)||_{C^{m-j,\alpha}}
+ ||h||_{C^{0}}).
\end{equation}
From \eqref{2.15}-\eqref{2.16} and the Bianchi identity, (as in \eqref{2.13}), one has
$\delta \hat L(h) = 2\delta \delta^{*}(\delta h)$ and the operator $\delta \delta^{*}$ is
elliptic with respect to Dirichlet boundary conditions. Since the boundary data $\delta h$
in \eqref{3.9} is included in the boundary operators $B_{j}$, this proves \eqref{3.11}.

  Using this and taking the trace of \eqref{3.11a} shows that
$$||D^{2}tr h||_{C^{m-2,\alpha}} \leq C(||\hat L(h)||_{C^{m-2,\alpha}} +
||B_{j}(h)||_{C^{m-j,\alpha}} + ||NN(tr h)||_{C^{m-2,\alpha}} + ||h||_{C^{0}}),$$
so that it suffices to prove that the boundary conditions $B$ cover $NN(tr h)$. For this,
a simple computation using \eqref{3.7}, (cf.~also \eqref{3.17} below),
gives
\begin{equation}\label{3.12a}
N(tr h) = 2H'_{h} - \delta((h(N))^{T}) - (\delta h)(N) + O(h),
\end{equation}
where $O(h)$ is of differential order 0 in $h$. Using the standard interpolation $||h||_{C^{m-1,\alpha}}
\leq \varepsilon ||h||_{C^{m,\alpha}}+ \varepsilon^{-1}||h||_{C^{0}}$, where $\varepsilon>0$ is an arbitrary constant, shows that it suffices
here and below only to consider terms with the leading number of derivatives of $h$.

  Now the Gauss equations at $\partial \mathcal{M}$ are $|A|^{2} - H^{2} + s_{\gamma_{\mathcal{M}}} = s_{g_{\mathcal{M}}}
- 2Ric(N,N)$ and hence,
$$(|A|^{2} - H^{2} + s_{\gamma_{\mathcal{M}}})'_{h} = -2\hat L(h)(N, N) + 2\delta^{*}\delta(h)(N, N) +
O(h).$$
One has $s_{\gamma_{\mathcal{M}}}'(h^{T}) = -\Delta tr h^{T} + \delta \delta (h^{T}) + O(h^{T})$
and $A'_{h}$, $H'_{h}$ only involve first order derivatives in $h$. Writing then $h^{T}
= B_{0}(h) + \frac{1}{n}tr_{\partial \mathcal{M}}h \gamma_{\mathcal{M}}$, it follows that $tr_{\partial \mathcal{M}}h$
at $\partial \mathcal{M}$ is controlled by $\hat L (h)$, $B_{j}(h)$, in that
$$||tr_{\partial \mathcal{M}}h||_{C^{m,\alpha}} \leq C(||h||_{C^{m-1,\alpha}} +
||\hat L(h)||_{C^{m-2,\alpha}} +  ||B_{j}(h)||_{C^{m-j,\alpha}}),$$
and hence
$$||h^{T}||_{C^{m,\alpha}} \leq C(||h||_{C^{m-1,\alpha}} +
||\hat L(h)||_{C^{m-2,\alpha}} +  ||B_{j}(h)||_{C^{m-j,\alpha}}),$$
i.e.~$h^{T}$ is controlled at $\partial \mathcal{M}$ by $\hat L(h)$ and $B_{j}(h)$. Next, at
$\partial \mathcal{M}$, one has $-(\delta h)(T) = \nabla_{N}h(N, T) + \nabla_{e_{i}}h(e_{i},T)$,
which then gives control as above on $(\nabla_{N}h)(N, T)$, and so control on
$\nabla_{N}(h(N)^{T})$. In turn, this gives then control on $\delta_{\partial \mathcal{M}}
(\nabla_{N}(h(N)^{T}))$, which modulo lower order (curvature) terms, equals
$N(\delta(h(N)^{T}))$.  The $N$-derivative of \eqref{3.12a} also holds and shows
that control of $N(\delta(h(N)^{T}))$ implies control of $NN(tr h)$, so that
\eqref{3.12} holds, provided $N(H'_{h})$ is controlled. But the Riccati equation gives
$N(H) = -|A|^{2} - Ric(N,N)$; taking the linearization of this in the direction $h$
shows that $N(H'_{h})$ is indeed controlled by $\hat L(h)$ and the boundary
conditions $B_{j}$. This completes the proof of ellipticity.

   Next, we prove the operator $\hat L$ with boundary conditions \eqref{3.9} is
self-adjoint. To begin, integrating the terms in the expression \eqref{2.16} for
$\hat L$ by parts over $\mathcal{M}$ gives
$$\int_{\mathcal{M}}\langle D^{*}D(h), k\rangle + \int_{\partial \mathcal{M}}\langle \nabla_{N}h, k\rangle =
\int_{\mathcal{M}}\langle D^{*}D(k), h\rangle + \int_{\partial \mathcal{M}}\langle \nabla_{N}k, h\rangle,$$
$$\int_{\mathcal{M}}\delta\delta h tr k + \int_{\partial \mathcal{M}}(\delta h)(N)trk =
\int_{\mathcal{M}}\langle h, D^{2}(trk)\rangle - \int_{\partial \mathcal{M}}h(N, dtrk),$$
and
$$\int_{\mathcal{M}}(\Delta trh) tr k - \int_{\partial \mathcal{M}}N(trh)trk =
\int_{\mathcal{M}}(\Delta trk)tr h - \int_{\partial \mathcal{M}}N(trk) tr h.$$
Here the boundary terms on $S(r)$ all tend to 0 as $r \rightarrow \infty$, since
the components of $h$ and $k$ are in $C_{\delta}^{m,\alpha}$ and $\delta >
\frac{1}{2}$. It follows that
\begin{equation}\label{3.13}
\int_{\mathcal{M}}\langle \hat L(h), k \rangle + \int_{\partial \mathcal{M}}\langle B(h), k \rangle =
\int_{\mathcal{M}}\langle \hat L(k), h \rangle + \int_{\partial \mathcal{M}}\langle B(k), h \rangle,
\end{equation}
where
\begin{equation}\label{3.14}
\langle B(k), h \rangle = \langle \nabla_{N}k, h \rangle + h(N, dtr k)
- (\delta k)(N) tr h - tr h N(tr k).
\end{equation}
Setting $Z(k, h) = \langle B(k), h \rangle - \langle B(h), k \rangle$, we thus
need to show that
\begin{equation}\label{3.15}
\int_{\partial \mathcal{M}}Z(h, k) = 0,
\end{equation}
when $h$, $k$ satisfy the boundary conditions \eqref{3.9}.

  Thus suppose $h$ and $k$ both satisfy \eqref{3.9}. A simple calculation shows
that $(\delta k)(T) = 0$ is equivalent to
\begin{equation}\label{3.16}
(\nabla_{N}k)(N)^{T} = \delta_{\partial \mathcal{M}}(k^{T}) - \alpha(k(N)),
\end{equation}
where $\alpha(k(N)) = [A(k(N)) + Hk(N)^{T}]$, (all taken on $\partial \mathcal{M}$),
while $(\delta k)(N) = 0$ is equivalent to
\begin{equation}\label{3.17}
N(k_{00}) = \delta_{\partial \mathcal{M}}(k(N)^{T}) + \langle A, k \rangle - k_{00}H.
\end{equation}
The same equations hold for $h$, and one also has
\begin{equation}\label{3.18}
h^{T} = \phi_{h}\gamma, \ \ {\rm and} \ \ k^{T} = \phi_{k}\gamma.
\end{equation}
We thus need to calculate
$$B(k, h) = \langle \nabla_{N}k, h \rangle + h(N, dtr k) - tr h N(tr k),$$
and skew-symmetrize. To begin, write $\langle \nabla_{N}k, h\rangle =
\langle (\nabla_{N}k)(N), h(N)\rangle + \langle (\nabla_{N}k)(e_{i}),
h(e_{i})\rangle$, so that $\langle \nabla_{N}k, h\rangle = N(k_{00})h_{00} +
\phi_{h}[N(tr k) - N(k_{00})] + 2\langle (\nabla_{N}k)(N), h(N)^{T} \rangle$, where
we have used the relation $tr_{\gamma}\nabla_{N}k = tr_{N}\nabla_{N}k - N(k_{00})$.
Thus, $B(k, h)$ equals
\begin{equation}\label{3.19}
N(k_{00})h_{00} + \phi_{h}[N(tr k) - N(k_{00})] + 2\langle (\nabla_{N}k)(N), h(N)^{T}
\rangle - N(tr k)[tr h - h_{00}] + \langle h(N)^{T}, dtr k \rangle.
\end{equation}
By \eqref{3.16} and \eqref{3.18},
$$2\langle (\nabla_{N}k)(N), h(N)^{T} \rangle = -2\langle d\phi_{k}, h(N)^{T}\rangle
- 2\alpha(k, h) = -2\phi_{k}\delta_{\partial \mathcal{M}}(h(N)^{T}) - 2\alpha(k, h) + \omega_{1}$$
where $\omega_{1}$ is a divergence term and $\alpha(k, h) = \langle \alpha(k(N)),
h(N)^{T} \rangle$. Similarly, by \eqref{3.17} and \eqref{3.18},
$$N(k_{00}) = \delta_{\partial \mathcal{M}}(k(N)^{T}) + (\phi_{k} - k_{00})H,$$
where here and in the following $\delta = \delta_{\partial \mathcal{M}}$. Note also that
$\langle h(N)^{T}, dtr k \rangle = tr k \delta_{\partial N}(h(N)^{T}) + \omega_{2}$,
where $\omega_{2}$ is another divergence term. Since \eqref{3.15} involves
integration over $\partial \mathcal{M}$, in the following we ignore the divergence terms.
Substituting these computations in \eqref{3.19} gives
$$\delta(k(N)^{T})[h_{00} - \phi_{h}] + \delta(h(N)^{T})[tr k - 2\phi_{k}] -
(n-1)\phi_{h}N(tr k) + H(\phi_{k} - k_{00})(h_{00} - \phi_{h}) - 2\alpha(h, k).$$
When skew-symmetrizing, the last two terms $H(\phi_{k} - k_{00})(h_{00} - \phi_{h}) -
2\alpha(h, k)$ cancel, while the first three terms combine to give
$$-(n-1)[\phi_{h}\delta(k(N)^{T}) - \phi_{k}\delta(h(N)^{T})] - (n-1)[N(tr k)\phi_{h} -
N(tr h)\phi_{k}],$$
or equivalently, (after dividing by $n-1$),
\begin{equation}\label{3.20}
-\phi_{h}[N(tr k) + \delta(k(N)^{T})] + \phi_{k}[N(tr h) + \delta(h(N)^{T})].
\end{equation}

  On the other hand, by \eqref{3.6} or \eqref{3.7},
$$2(H')_{k} = tr[\nabla_{N}k - 2\delta^{*}(k(N)^{T}) - \delta^{*}(k_{00}N)]$$
$$= N(tr k) + 2\delta(k(N)^{T}) - k_{00}H - N(k_{00}).$$
Substituting \eqref{3.17} gives
$$2(H')_{k} = N(tr k) + \delta(k(N)^{T}) - H\phi_{k},$$
so \eqref{3.20} becomes
$$-\phi_{h}[2(H')_{k} + H\phi_{k}] + \phi_{k}[2(H')_{h} + H\phi_{h}] = -2\phi_{h}(H')_{k} +
2\phi_{k}(H')_{h}.$$
This vanishes exactly when $H_{k}'$ and $H_{h}'$ vanish. This completes the proof.

{\endproof}

 The main step in the proof of the manifold theorem, (Theorem 3.6), is the following result.
\begin{theorem} \label{t3.3}
Suppose $\pi_{1}(M, \partial M) = 0$ and $m \geq 3$. Then at any $\widetilde g \in
{\mathbb E}_{S}$, the map $\hat \Phi$ is a submersion, i.e.~the derivative
\begin{equation} \label{3.21}
(D\hat \Phi)_{\widetilde g}: T_{\widetilde g}Met_{D}^{m,\alpha}(\mathcal{M}) \rightarrow
T_{\hat \Phi (\widetilde g)}S_{\delta}^{m-2,\alpha}(\mathcal{M})
\end{equation}
is surjective and its kernel splits in $T_{\widetilde g}Met_{D}^{m,\alpha}(\mathcal{M})$.
\end{theorem}

{\bf Proof:}
  The operator $\hat L = 2D\hat \Phi_{\widetilde g}$ is elliptic in the interior,
and the boundary data in Lemma 3.2 give a self-adjoint elliptic boundary value problem.
Let $S_{B}^{m,\alpha}(\mathcal{M})$ be the space of $C^{m,\alpha}$ symmetric bilinear forms on
$\mathcal{M}$ satisfying the boundary condition $B(h) = 0$ from Lemma 3.2, i.e.
$$B(h) = \{\delta h, [h^{T}]_{0}, (H')_{h}\} = (0,0,0).$$
Clearly, $S_{B}^{m,\alpha}(N) \subset S_{D}^{m,\alpha}(N)$, where
$S_{D}^{m,\alpha}(\mathcal{M}) = T_{\widetilde g}(Met_{D}^{m,\alpha}(\mathcal{M}))$. Throughout the
following, we set $\widetilde g = g$. The operator $\hat L$, mapping
$$S_{B}^{m,\alpha}(\mathcal{M}) \rightarrow S_{\delta}^{m-2,\alpha}(\mathcal{M}),$$
$$\hat L(h) = f, \ \ B(h) = 0 \ \ {\rm at} \ \ \partial \mathcal{M},$$
is then Fredholm, of Fredholm index 0. On $S_{B}^{m,\alpha}(\mathcal{M})$, the image
$Im (\hat L)$ is a closed subspace of the range $S^{m-2,\alpha}(\mathcal{M})$, of finite
codimension, and with codimension equal to dimension of the kernel $K$.

  If $K = 0$, then $\hat L$ maps $S_{B}^{m,\alpha}(\mathcal{M})$ onto $S_{\delta}^{m-2,\alpha}(\mathcal{M})$,
which proves the result. Thus suppose $K \neq 0$. Then as in \eqref{3.13},
by the self-adjointness, one has for any $h \in S_{B}^{m,\alpha}(\mathcal{M})$ and $k \in K$,
$$\int_{\mathcal{M}}\langle \hat L(h), k \rangle = \int_{\mathcal{M}}\langle h, \hat L(k) \rangle = 0,$$
since the boundary terms vanish and $\hat L(k) = 0$. Thus $Im(\hat L|_{S_{B}^{m,\alpha}(\mathcal{M})})
= K^{\perp}$, (where $K^{\perp}$ is taken with respect to the $L^{2}$ inner
product). To prove surjectivity on $S_{D}^{m,\alpha}(\mathcal{M})$, it thus suffices to prove
that for any $k \in K$, there exists $h \in S_{D}^{m,\alpha}(\mathcal{M})$ such that
\begin{equation}\label{3.22}
\int_{\mathcal{M}}\langle \hat L(h), k \rangle \neq 0.
\end{equation}

  Suppose then \eqref{3.22} does not hold, so that
\begin{equation}\label{3.23}
\int_{\mathcal{M}}\langle \hat L(h), k \rangle = 0,
\end{equation}
for all $h \in S_{D}^{m,\alpha}(\mathcal{M})$, i.e.~for which $\delta h = 0$ on $\partial \mathcal{M}$.
Integrating by parts, it follows that
\begin{equation}\label{3.24}
\int_{\mathcal{M}}\langle h, \hat L(k) \rangle + \int_{\partial \mathcal{M}}Z(h, k) = 0,
\end{equation}
for $Z(h, k)$ as following \eqref{3.14}. As before, the boundary terms at
infinity vanish, since $\delta > \frac{1}{2}$.

  Choosing $h \in S_{D}^{m,\alpha}(\mathcal{M})$ arbitrary of compact support in $\mathcal{M}$,
it follows from \eqref{3.24} that
\begin{equation} \label{3.25}
\hat L(k) = 0,
\end{equation}
i.e.~$k$ is formally tangent to $\hat Z = \hat \Phi^{-1}(0)$. Of course this is
already known, since $k \in K$. Moreover, one also has
\begin{equation} \label{3.26}
\delta k = 0 \ \ {\rm on} \ \ \mathcal{M}.
\end{equation}
To see this, let $h = \delta^{*}V$, with $V$ any vector field vanishing on $\partial \mathcal{M}$.
Since $g$ is Einstein and so $(Ric - \frac{s}{2}g)'_{\delta^{*}V} = 0$, it follows
from \eqref{3.5} and \eqref{3.6} that $\hat L(h) = \delta^{*}Y$, where $Y =
2\delta \delta^{*}V$. As in Lemma 2.2, the operator $\delta \delta^{*}$ is surjective,
(in fact an isomorphism), on vector fields vanishing at $\partial \mathcal{M}$, so that $Y$ may be
arbitrarily prescribed. Moreover, $h \in S_{D}^{m,\alpha}(\mathcal{M})$ if and only if $Y = 0$ at
$\partial \mathcal{M}$. Then \eqref{3.23} gives
$$0 = \int_{\mathcal{M}}\langle \hat L(\delta^{*}V), k \rangle = \int_{\mathcal{M}}\langle \delta^{*}Y, k \rangle
= \int_{\mathcal{M}}\langle Y, \delta k \rangle + \int_{\partial \mathcal{M}}k(Y, N) = \int_{\mathcal{M}}\langle Y,
\delta k \rangle,$$
since $Y = 0$ on $\partial \mathcal{M}$. Here we have used again the fact that the boundary term
at infinity vanishes, since $|k| = O(r^{-\delta})$ and $|Y| = O(r^{-1-\delta})$. Since
$Y$ is otherwise arbitrary, this gives \eqref{3.26}.

  Returning now to \eqref{3.24}, \eqref{3.25} gives
\begin{equation} \label{3.27}
\int_{\partial \mathcal{M}}Z(h, k) = 0,
\end{equation}
for all $h$ with $\delta h = 0$ on $\partial \mathcal{M}$. Next, we choose certain test forms
$h \in S_{D}(\mathcal{M})$ in \eqref{3.27}. First, choose $h$ such that $h = 0$ on $\partial \mathcal{M}$.
Then $\nabla_{N}h$ is freely specifiable, subject to the divergence constraint $\delta h = 0$;
all computations here and below are at $\partial \mathcal{M}$. Since $h = 0$, this constraint gives
$(\nabla_{N}h)(N) = 0$, which is equivalent to the tangential and normal constraints:
\begin{equation} \label{3.28}
(\nabla_{N}h)(N, T) = 0,
\end{equation}
\begin{equation} \label{3.29}
N(h_{00}) = 0,
\end{equation}
for any $T$ tangent to $\partial \mathcal{M}$. Choosing $h$ and $\nabla_{N}h$
satisfying $h = 0$ and \eqref{3.28}-\eqref{3.29} at $\partial \mathcal{M}$, the terms
$(\nabla_{N}h)(T_{1}, T_{2})$ are freely specifiable on $\partial \mathcal{M}$,
where $T_{1}, T_{2}$ are any vectors tangent to $\partial \mathcal{M}$.
Substituting such $h$ in \eqref{3.27} and using \eqref{3.26}, it follows that
\begin{equation} \label{3.30}
\int_{\partial \mathcal{M}}\langle \nabla_{N}h, k\rangle  + (k_{00}-trk)N(trh) = 0.
\end{equation}
Now choose $\nabla_{N}h = fg^{T},$ where $g^{T} = g|_{T(\partial \mathcal{M})}$.
This choice satisfies the constraints \eqref{3.28}-\eqref{3.29}. The integrand
in \eqref{3.30} then becomes $ftr^{T}k - N(trh)tr^{T}k$. Since $N(trh) =
\langle \nabla_{N}h, g\rangle = nf$, and since $f$ is arbitrary, it
follows that $tr^{T}k = 0$. In turn, since the tangential part of
$\nabla_{N}h$ is arbitrary, \eqref{3.30} implies
\begin{equation} \label{3.31}
k^{T} = 0, \ \ {\rm on} \ \ \partial \mathcal{M}.
\end{equation}

\begin{lemma}\label{l3.4}
At $\partial \mathcal{M},$ one has
\begin{equation} \label{3.32}
(A'_{k})^{T} = 0,
\end{equation}
i.e.~ $(\nabla_{N}k)^{T} = 2[\delta^{*}(k(N)^{T})]^{T} + k_{00}A$, since
$k^{T} = 0$, cf.~\eqref{3.6}.
\end{lemma}

{\bf Proof:} The proof is a straightforward, but rather long computation.
To begin, as preceding \eqref{3.19} and using \eqref{3.31}, one has
$\langle \nabla_{N}h, k\rangle = 2\langle (\nabla_{N}h)(N), k(N)^{T}\rangle
+ N(h_{00})k_{00}$. By \eqref{3.16}, $(\nabla_{N}h)(N)^{T} = \delta_{\partial \mathcal{M}}(h^{T})
- \alpha (h(N))$, so that
\begin{equation}\label{e3.33}
\int_{\partial \mathcal{M}}\langle \nabla_{N}h, k\rangle  = \int_{\partial \mathcal{M}}2\langle
\delta_{\partial \mathcal{M}}(h^{T}), k(N)^{T}\rangle  + N(h_{00})k_{00} - 2\alpha (h, k)
\end{equation}
$$= \int_{\partial \mathcal{M}}2\langle h^{T}, (\delta_{\partial \mathcal{M}})^{*}(k(N)^{T})\rangle  +
N(h_{00})k_{00} - 2\alpha(h, k).$$
Further, for $Z$ tangent to $\partial \mathcal{M}$, one has
$(\delta_{\partial \mathcal{M}})^{*}(k(N)^{T})(Z,Z) = \langle \nabla^{T}_{Z}k(N)^{T}, Z\rangle  =
\langle \nabla_{Z}k(N)^{T}, Z\rangle  = \delta^{*}(k(N)^{T})(Z,Z)$, where now
$\delta^{*}$ is taken with respect to the ambient metric $g_{\mathcal{M}}$, (not the
boundary metric $\gamma_{\mathcal{M}}$). So this gives
\begin{equation} \label{e3.34}
\int_{\partial \mathcal{M}}\langle \nabla_{N}h, k\rangle  = \int_{\partial \mathcal{M}}\langle h^{T},
2\delta^{*}(k(N)^{T})\rangle + N(h_{00})k_{00} - 2\alpha(h, k).
\end{equation}

 On the other hand, one computes $\langle \nabla_{N}k, h\rangle  = \langle (\nabla_{N}k)^{T},
h^{T}\rangle + \langle \nabla_{N}k(N), h(N)\rangle  = \langle (\nabla_{N}k)^{T}, h^{T}\rangle -
\langle \alpha (k(N)), h(N)^{T}\rangle + N(k_{00})h_{00}$, again by \eqref{3.16} and \eqref{3.31}.
Taking the difference of this with \eqref{e3.34} and noting that $\alpha$ is symmetric, gives
\begin{equation} \label{e3.35}
\int_{\partial \mathcal{M}}\langle h^{T}, (\nabla_{N}k)^{T} - 2\delta^{*}(k(N)^{T})\rangle  +
N(k_{00})h_{00} - N(h_{00})k_{00}  = E,
\end{equation}
where via \eqref{3.14}-\eqref{3.15}, $E$ is given by
$$E = \int_{\partial \mathcal{M}}[k(N,dtrh) - h(N,dtrk)] - [N(trh)trk - trh N(trk)]. $$
Computing this term-by-term gives: $k_{00}N(trh) + \langle k(N)^{T},
d^{T}trh\rangle - h_{00}N(trk) - \langle h(N)^{T},d^{T}trk\rangle - N(trh)trk +
trh N(trk)$. Since $trk = k_{00}$, the first and second-to-last terms
cancel. Integrating over $\partial \mathcal{M}$ and using the divergence theorem
shows that
\begin{equation} \label{e3.36}
E = \int_{\partial \mathcal{M}}trh\delta_{\partial \mathcal{M}}(k(N)^{T}) -
k_{00}\delta^{T}(h(N)^{T}) - h_{00}N(trk) + trh N(trk).
\end{equation}
Next we claim that
\begin{equation} \label{e3.37}
\delta_{\partial \mathcal{M}}(h(N)^{T}) = N(h_{00}) + Hh_{00} - \langle A, h\rangle ,
\end{equation}
and similarly for $k$. This follows from the following computation:
$\delta_{\partial \mathcal{M}}(h(N)^{T}) = \delta_{\partial \mathcal{M}}(h(N)) - \delta^{T}(h_{00}N) =
\delta_{\partial \mathcal{M}}(h(N)) + Hh_{00}$, while $\delta_{\partial \mathcal{M}}(h(N)) = \delta (h(N)) +
N(h_{00})$. Since $\delta (h(N)) = (\delta h)(N) - \langle A, h\rangle$,
this gives the claim. Substituting \eqref{e3.37} into \eqref{e3.36}, and using
$\langle A, k\rangle = 0$ implies that
$$E = \int_{\partial \mathcal{M}}trh(N(k_{00}) + Hk_{00}) - k_{00}(N(h_{00}) +
Hh_{00} - \langle A,h\rangle ) - h_{00}N(trk) + trh N(trk), $$
and rearranging terms gives
\begin{equation} \label{e3.38}
E = \int_{\partial \mathcal{M}}\langle A, h\rangle k_{00} + N(trk)[trh - h_{00}] +
trhN(k_{00}) - k_{00}N(h_{00}) + H(trhk_{00} - trkh_{00}).
\end{equation}
Now substitute \eqref{e3.38} into \eqref{e3.35}: the $k_{00}N(h_{00})$ term cancels
to give
\begin{equation}\label{e3.39}
\int_{\partial \mathcal{M}}\langle h^{T},(\nabla_{N}k)^{T}-2\delta^{*}(k(N))^{T}\rangle  -
\langle A, h\rangle k_{00} =
\end{equation}
$$-\int_{\partial \mathcal{M}}N(k_{00})h_{00} - N(trk)[trh - h_{00}] -
trhN(k_{00}) - H(trhk_{00} - trkh_{00}). $$
The integrand on the right combines to: $-N(k_{00})(h_{00} - trh) -
N(trk)[h_{00} - trh] - Htrk(h_{00} - trh) = -[N(k_{00}) + N(trk) +
Htrk](h_{00}-trh)$. Since $h_{00} - trh = -tr^{T}h = -\langle h^{T},
g^{T}\rangle$, and since $h^{T}$ may be chosen arbitrarily, (the constraint
$\delta h = 0$ imposes no constraint on $h^{T}$), it follows that
\begin{equation} \label{e3.40}
(\nabla_{N}k)^{T} = 2[\delta^{*}(k(N)^{T})]^{T}+ k_{00}A + [N(k_{00}) + N(trk)
+ Htrk]g^{T}.
\end{equation}

 To complete the proof of \eqref{3.32}, we thus need to show that
\begin{equation} \label{e3.41}
N(k_{00}) + N(trk) + Htrk = 0.
\end{equation}
To obtain \eqref{e3.41}, take the $g^{T}$-trace of \eqref{e3.40}. One has
$\langle \nabla_{N}k, g^{T}\rangle  = N(trk) - N(k_{00})$, while
$\langle \delta^{*}(k(N)^{T}), g^{T}\rangle  = \langle \nabla_{e_{i}}k(N)^{T},
e_{i}\rangle  = \langle \nabla_{e_{i}}k(N), e_{i}\rangle - k_{00}H =
\langle (\nabla_{e_{i}}k)(N), e_{i}\rangle  - k(\nabla_{e_{i}}N, e_{i}) -
k_{00}H = -N(k_{00}) - k_{00}H$, the last equality using \eqref{3.31} and
\eqref{3.26}. This gives
$$N(trk) - N(k_{00}) = -2N(k_{00}) - 2k_{00}H + k_{00}H - n[N(k_{00}) +
N(trk) + Htrk], $$
which implies \eqref{e3.41}. This completes the proof of the Lemma.

{\endproof}

  To complete the proof of Theorem 3.3, \eqref{3.31} and \eqref{3.32} show that
$$k^{T} = (A'_{k})^{T} = 0,$$
at $\partial \mathcal{M}$. One also has $\hat L(k) = \delta k = 0$ on $\mathcal{M}$, so that $k$ is
an infinitesimal Einstein deformation on $\mathcal{M}$. By the local unique continuation
result of \cite{AH}, together with the global hypothesis $\pi_{1}(M, \partial M)
= 0$, it follows that $k = 0$. This shows that $\hat L$ is surjective. The fact
that its kernel splits is standard, cf.~\cite{A3}. This completes the proof.

{\endproof}

  Via the implicit function theorem, one obtains:
\begin{corollary}\label{c3.5}
Suppose $\pi_{1}(M, \partial M) = 0$ and $m \geq 3$. Then the local spaces
${\mathbb E}_{D}^{m,\alpha}$ are infinite dimensional $C^{\infty}$ Banach manifolds,
with
\begin{equation} \label{e3.42}
T_{\widetilde g}{\mathbb E}_{D} = Ker(D\hat \Phi_{\widetilde g}).
\end{equation}
\end{corollary}

{\bf Proof:}
 This is an immediate consequence of Theorem 3.3, the fact from Proposition 2.1 that
${\mathbb E}_{D} = Z_{D}$, (cf.~\eqref{2.18}), and the implicit function theorem,
(or regular value theorem), in Banach spaces.
{\endproof}

  This leads to the main result of this section.

\begin{theorem}\label{t3.6}
Suppose $\pi_{1}(M, \partial M) = 0$ and $m \geq 3$. Then the moduli space
${\mathcal E}_{S} = {\mathcal E}_{S}^{m,\alpha}$ is a $C^{\infty}$ smooth
infinite dimensional Banach manifold for which the boundary map
\begin{equation}\label{3.43}
\Pi_{B}: {\mathcal E}_{S} \rightarrow Met^{m,\alpha}(\partial M)\times C^{m-1,\alpha}
(\partial M),
\end{equation}
is a $C^{\infty}$ smooth Fredholm map, of Fredholm index 0.
\end{theorem}

{\bf Proof:}
Recall from \S 1 that the moduli space ${\mathcal E}_{S}$ of static vacuum Einstein metrics
is defined to be the quotient ${\mathbb E}_{S}^{m,\alpha}/{\mathcal D}_{1}^{m+1,\alpha}$.
The local spaces ${\mathbb E}_{D}$ are smooth Banach manifolds and depend smoothly
on the background metric $\widetilde g$, since the divergence-free gauge condition
\eqref{2.18} varies smoothly with $\widetilde g$. As noted preceding Lemma 2.2, the
action of ${\mathcal D}_{1}$ on ${\mathbb E}$ is free and the local spaces ${\mathbb E}_{D}$
are smooth local slices for the action of ${\mathcal D}_{1}$ on  ${\mathbb E}_{S}$.
Hence the global space ${\mathbb E}_{S}$ is a smooth Banach manifold, as is the
quotient ${\mathcal E}_{S}$. The local slices ${\mathbb E}_{D}$ represent local
coordinate patches for ${\mathcal E}_{S}$.

  Proposition 3.1 implies that the boundary map $\Pi_{B}: {\mathbb E}_{S}^{m,\alpha}
\rightarrow Met^{m,\alpha}(\partial M)\times C^{m-1,\alpha}(\partial M)$ is smooth
and Fredholm, of Fredholm index 0. Moreover, $\Pi_{B}$ is invariant under the action of
${\mathcal D}_{1}^{m+1,\alpha}(M)$ on ${\mathcal E}_{S}^{m,\alpha}$ and so it descends
to a smooth Fredholm map as in \eqref{3.43}, still of index 0.

{\endproof}

  The boundary conditions $(\gamma, H)$ for the operator $\hat \Phi$ are also self-adjoint; in fact 
they arise naturally from a variational principle (Lagrangian) on a space of static (non-vacuum) metrics. 

   To describe this, let $S(h) = -\Delta tr h + \delta \delta(h) - \langle Ric,
h \rangle$ be the linearization of the scalar curvature $s$, with adjoint $S^{*}$
given by
$$S^{*}u = D^{2}u -\Delta u \, g - uRic.$$
It is well-known that the static vacuum equations are given by $S^{*}u = 0$ and
$s = 0$ on $M$.

\begin{proposition}\label{p3.7}
For $(M, g, u)$ as above, the Bartnik boundary conditions $(\gamma, H)$ give
a well-defined variational problem for the Lagrangian
\begin{equation} \label{3.44}
{\mathcal L}(g, u) = \int_{M}usdV_{g} - 16\pi m_{ADM}: Met_{\delta}^{m,\alpha}(M)
\times C_{\delta}^{m,\alpha}(M) \rightarrow {\mathbb R},
\end{equation}
where $m_{ADM}$ is the ADM mass of $(M, g, u)$. The gradient $\nabla {\mathcal L}$
of ${\mathcal L}$ at $(g, u)$ is given by
\begin{equation} \label{3.45}
\nabla {\mathcal L} = (S^{*}u + {\tfrac{1}{2}}usg, \ s, \ uA - N(u)\gamma,\  2u),
\end{equation}
in the following sense: if $(h, u')$ is a variation of $(g, u)$ inducing the
variation $(h^{T}, H'_{h})$ of the boundary data, then
\begin{equation} \label{3.46}
d{\mathcal L}(h, u', h^{T}, H'_{h}) = \int_{M}[\langle S^{*}u + {\tfrac{1}{2}}usg,
h \rangle + su'] + \int_{\partial M}[\langle uA - N(u) \gamma, h^{T} \rangle + 2uH'_{h}].
\end{equation}
   In particular, the static vacuum equations are critical points for
${\mathcal L}$ with data $(\gamma, H)$ fixed.
\end{proposition}

{\bf Proof:} Suppose that $D$ is compact domain in $M$, with $N$ the outward unit
normal from $D$. Varying $(g, u)$ in the direction $(h, u')$ then gives
\begin{equation} \label{3.47}
D{\mathcal L}_{0}(h, u') = \int_{D}(us' + u's + us(dV)') = \int_{D}\langle u, s(h) \rangle
+ {\tfrac{1}{2}}us \langle g, h \rangle + su',
\end{equation}
where ${\mathcal L}_{0} = \int_{D}usdV_{g}$. A straightforward integration by parts gives
\begin{equation} \label{3.48}
\int_{D}\langle S(h), u \rangle = \int_{D}\langle h, S^{*}u \rangle +
\int_{\partial D}-uN(tr h) - (\delta h)(N)u - \langle h(N), du \rangle + tr h N(u).
\end{equation}

  The equations \eqref{3.47} and \eqref{3.48} imply immediately the bulk Euler-Lagrange
equations - the first two terms in \eqref{3.45}. If the bulk term (over $D$) vanishes,
then since $u'$ is arbitrary $s = 0$, so this gives
$$S^{*}u = 0,$$
with $s = 0$, which are the static vacuum equations.

  For the boundary terms, from \eqref{3.7} one has
$$2H'_{h} = N(tr h) + 2\delta(h(N)^{T}) - h_{00}H - N(h_{00}).$$
Also by a simple calculation
$$(\delta h)(N) = \delta(h(N)^{T}) + \langle A, h \rangle - h_{00}H - N(h_{00}),$$
so that,
$$2H'_{h} - (\delta h)(N) = N(tr h) + \delta(h(N)^{T}) - \langle A, h \rangle.$$
This gives
$$\int_{\partial D}u(-N(tr h) - (\delta h)(N)) = \int_{\partial D}-2uH'_{h} +
\langle du, h(N)^{T}\rangle - u\langle A, h \rangle.$$
It follows that the boundary term in \eqref{3.48} is given by
\begin{equation} \label{3.49}
\int_{\partial D}-2uH'_{h} - u\langle A, h \rangle - N(u)h_{00} + N(u)tr h =
\int_{\partial D}-2uH'_{h} - u\langle A, h^{T} \rangle +
N(u)\langle h^{T}, \gamma \rangle.
\end{equation}

  Now let the outer boundary of $D$ equal $S(r)$ and consider the limit $r \rightarrow \infty$.
Then $u \rightarrow 1$ and $\int_{S(r)}N(u)\langle h^{T}, \gamma \rangle \rightarrow 0$.
It follows that
$$\lim_{r \rightarrow \infty}\int_{S(r)}-2uH'_{h} - u\langle A, h^{T} \rangle =
\lim_{r \rightarrow \infty}\int_{S(r)}(-N(tr h) - (\delta h)(N)) = 16\pi (m_{ADM})',$$
where the second equality follows from standard formulas for the ADM mass and its
variation, cf.~\cite{RT}, \cite{B3}. Here the variation $(m_{ADM})'$ is taken
in the direction $h$. The formula \eqref{3.48}-\eqref{3.49} is also valid at the
inner boundary $\partial M$, with respect to the outward normal. Changing to the inner
normal changes the sign of each term, and \eqref{3.45} and \eqref{3.46} then follow
immediately.

{\endproof}

  On-shell, i.e.~on the space of solutions ${\mathcal E}$, the Lagrangian
$${\mathcal L} = -16\pi m_{ADM} : {\mathcal E} \rightarrow {\mathbb R}$$
is a smooth function whose derivative is given by the boundary term in \eqref{3.46},
a result basically due to Bartnik \cite{B3}. After writing this work, we learned that
a special case of Proposition 3.7 has been noted, without proof, in a paper of
Miao, cf.~\cite{M2}.

\section{Curvature estimates and Properness of $\Pi^{o}$.}
\setcounter{equation}{0}

   In contrast to the previous section, where most of the computations were done on
the 4-manifold $(\mathcal{M}, g_{\mathcal{M}})$, in this section we focus mostly on the 3-dimensional
data $(M, g, u)$. Let $inj_{\partial M}$ denote the injectivity radius of the normal
exponential map from $\partial M$ in $M$ and let $R$ denote the (full) curvature
tensor of $g = g_{M}$. The main result of this section is the following collection of
apriori estimates.

\begin{theorem}\label{t5.1}
For $(M, g, u) \in {\mathcal E}^{o} = ({\mathcal E}^{m,\alpha})^{o}$ with $m \geq 2$,
one has a global pointwise estimate
\begin{equation}\label{5.1}
|R| \leq \Lambda,
\end{equation}
on $M$, where $\Lambda$ depends only on bounds for the boundary data $(\gamma, H)$.
Moreover, at $\partial M$, one has the bounds
\begin{equation}\label{5.1a}
|A| \leq \Lambda, \ \ inj_{\partial M} \geq \Lambda^{-1}.
\end{equation}
The estimates \eqref{5.1}, \eqref{5.1a} also hold for higher derivatives of $R$ and $A$,
up to order $m-2$, $m-1$ respectively.
\end{theorem}

{\bf Proof:} For points in the interior of $M$, of bounded distance away from $\partial M$,
this follows directly from the apriori interior estimates in \cite{A1} which state
\begin{equation}\label{5.2}
|R|(x) \leq \frac{K}{t^{2}(x)}, \ \ |d\log u|(x) \leq \frac{K}{t(x)},
\end{equation}
where $t(x) = dist(x, \partial M)$, where $K$ is an absolute constant. Similar
(scale-invariant) estimates hold for all higher derivatives of $R$ and $\log u$. So one
only needs to consider the behavior near $\partial M$. At $\partial M$, the Gauss and
Gauss-Codazzi (constraint) equations are given by:
\begin{equation}\label{5.3}
|A|^{2} - H^{2} + s_{\gamma} = -2R_{NN},
\end{equation}
\begin{equation}\label{5.4}
\delta(A - H\gamma) = -u^{-1}D^{2}u(N, \cdot).
\end{equation}
Also, $-R_{NN} = -Ric(N,N) = -u^{-1}NN(u) = u^{-1}(\Delta_{\partial M} u + HN(u))$,
so that
\begin{equation}\label{5.5}
u(|A|^{2} - H^{2} + s_{\gamma}) = 2(\Delta_{\partial M} u + HN(u)).
\end{equation}

  From \eqref{5.3}, a bound on $|R|$ at $\partial M$ gives immediately a bound on
$|A|$ at $\partial M$, given control of $(\gamma, H)$. Similarly, a bound on $|R|$
on $M$ gives a lower bound on the distance $d_{con}$ to the conjugacy locus of the
normal exponential map $exp_{\partial M}$.

  Now, again under a bound on $|R|$, the outer-minimizing property \eqref{1.10} implies
a lower bound on the distance $\delta_{\partial M}$ to the cut locus of $exp_{\partial M}$.
To see this, suppose that $\delta_{\partial M} << 1$ but $\Lambda$ in \eqref{5.1} is
bounded, $\Lambda \sim 1$. Then since $d_{con}$ is bounded below, there is a
geodesic $\zeta$ of length $2\delta_{\partial M}$ in $M$ meeting $\partial M$
orthogonally at points $p_{1}, p_{2}$. Let $T$ be the boundary of the tubular
neighborhood of $\zeta$ of radius $r$. Then $T$ intersects $\partial M$ in
the boundary of two discs $D_{1}$, $D_{2}$ of radius approximately $r$, (for $r$ small).
If $\delta_{\partial M} << r$, then $area T < area(D_{1}\cup D_{2})$. Further,
$T \subset M$ is homologous to $D_{1}\cup D_{2}$ in $M$. Removing then $D_{1}\cup D_{2}$
from $\partial M$ and attaching $T$ shows that $\partial M$ is not outer-minimizing in
$M$, giving a contradiction.

   Thus it suffices to obtain a curvature bound at or arbitrarily near $\partial M$.
The higher derivative estimates may then be obtained by standard elliptic regularity methods.
The proof of \eqref{5.1} is by a blow-up argument. If the curvature bound in \eqref{5.1} is
false, then there is a sequence $(M, g_{i}, u_{i}, x_{i}) \in {\mathcal E}^{o}$ with
bounded Bartnik boundary data such that
$$|R_{g_{i}}|(x_{i}) \rightarrow \infty.$$
Without loss of generality, assume that the curvature of $g_{i}$ is maximal at $x_{i}$.
We then rescale the metrics $g_{i}$ to $g_{i}'$ so that $|R|$ is bounded, and equals 1
at $x_{i}$,
\begin{equation}\label{5.5a}
|R_{g_{i}'}|(x_{i}) = 1, \ \ |R_{g_{i}'}|(y_{i}) \leq 1,
\end{equation}
for any $y_{i} \in (M, g_{i}')$. Thus define $g_{i}' = \lambda_{i}^{2}g_{i}$ where $\lambda_{i} = 
|R_{g_{i}}|(x_{i})$. This gives \eqref{5.5a} as well as $H_{g_{i}'} = \lambda_{i}^{-1}H_{g_{i}}$, 
$s_{\gamma_{i}'} = \lambda_{i}^{-2}s_{\gamma_{i}}$ and $t_{i}' = \lambda_{i}t_{i}$. Note that 
\eqref{5.2} implies that $t_{i}'(x_{i}) \leq \sqrt{K}$ so that $x_{i}$ remains within a uniformly bounded 
distance to the boundary $\partial M$ with respect to $g_{i}'$.  

   One may also need to rescale the potential $u$. For
reasons that will be clearer below, choose points $y_{i} \in M$ such that
$dist_{g_{i}'}(y_{i}, \partial M) = 1$ and $dist_{g_{i}'}(y_{i}, x_{i}) \leq  \sqrt{K}$,
and rescale $u_{i}$ so
\begin{equation}\label{5.6}
u_{i}'(y_{i}) = 1.
\end{equation}

  The sequence $(M, g_{i}', u_{i}')$ has uniformly bounded curvature and uniform control
of the boundary geometry, (boundary metric, $2^{\rm nd}$ fundamental form and normal
exponential map). By \eqref{5.6} and the Harnack inequality for positive harmonic functions, 
the potential $u_{i}'$ is also uniformly bounded in domains of bounded diameter about $y_{i}$ or 
$x_{i}$. It follows from the convergence theorem in \cite{AT} for manifolds-with-boundary that a 
subsequence converges weakly, (i.e.~in $C^{1,\alpha}$), to a $C^{1,\alpha}$ static limit $(X, g, u, x)$ 
with boundary $(\partial X, \gamma, u)$. The convergence is uniform on compact subsets. More precisely, 
given any smooth compact domain $\Omega$ in the manifold-with-boundary $X$, there is a subsequence, 
also denoted $\{i\}$ and embeddings $F_{i}: \Omega \rightarrow (M_{i}, g_{i}', u_{i}')$ such that 
$F_{i}^{*}(g_{i}', u_{i}') \rightarrow (\Omega, g, u)$ in the $C^{1,\alpha}$ topology. Moreover, for 
$\Omega_{1} \subset \Omega_{2}$, $g_{2}|_{\Omega_{1}} = g_{1}$. For the proof, we refer to 
[AT, Theorem 3.1], and more precisely to the local or pointed version of this result in [AT, Theorem 3.1.1]. 

   By the normalization in \eqref{5.5a}, the limit $(X, g)$ is complete (without singularities) up to the 
boundary $\partial X$. Since $\partial M$ is outer-minimizing in $(M, g_{i})$, the $C^{0}$ convergence to the limit
implies that $\partial X$ is weakly outer-minimizing in $X$: if $D$ is any compact smooth domain
in $\partial X$ and $D'$ is a surface in $X$ with $\partial D' = \partial D$, then
\begin{equation}\label{5.7}
area D' \geq area D.
\end{equation}

   One has $\partial X = {\mathbb R}^{2}$, the boundary metric $\gamma$ is flat, $H = 0$,
so $\partial X$ is a minimal surface in $X$. One has $u > 0$ in the interior of $X$,
(by the maximum principle), but may have $u = 0$ somewhere or everywhere on
$\partial X$. The bound \eqref{5.5a} and the static equations imply that $u_{i}$ is
uniformly bounded up to $\partial X$, within bounded distance to $x_{i}$ and the
limit potential $u$ extends at least $C^{1,\alpha}$ up to $\partial X$.

  We will prove below that the convergence to the limit is smooth, so that in particular
\begin{equation}\label{5.8}
|R|(x) = 1,
\end{equation}
where $x = \lim x_{i}$ and $R = R_{X}$.

  On the blow-up limit $(X, g)$, \eqref{5.5} holds and becomes
$${\tfrac{1}{2}}u|A|^{2} = \Delta_{\partial M}u,$$
on $\partial X$. This equation holds weakly on $\partial X$ with $u \in C^{1,\alpha}
(\partial X)$; elliptic regularity then implies it holds strongly, and
$u \in C^{3,\alpha}(\partial X)$. Since $u$ is harmonic, $u$ is thus $C^{3,\alpha}$
up to $\partial X$. Also, by the Riccati equation $N(H) = -|A|^{2} -
R_{NN} = -|A|^{2} + \frac{1}{2}(|A|^{2} - H^{2} + s_{\gamma})$, so that
\begin{equation}\label{5.9}
N(H) = -{\tfrac{1}{2}}(|A|^{2} + H^{2} - s_{\gamma}).
\end{equation}
This holds pointwise on the blow-up sequence $(M, g_{i}', u_{i})$ and since
$s_{\gamma} \rightarrow 0$ and $H \rightarrow 0$ for $g_{i}'$, it follows that $N(H)$
is defined pointwise on the limit $\partial X$ and on $\partial X$,
\begin{equation}\label{5.10}
N(H) \leq 0,
\end{equation}
with equality on any domain only when $A = 0$.

   Since $\partial X$ is minimal, \eqref{5.10} and the outer-minimizing property
\eqref{5.7} imply that
\begin{equation}\label{5.11}
N(H) = 0,
\end{equation}
on $\partial X$. In more detail, \eqref{5.7} and the fact that $H = 0$ on $\partial X$
implies the $2^{\rm nd}$ order stability of $\partial X$, in that the $2^{\rm nd}$
variation of the area of $\partial X$ is non-negative. Thus, for all $f$ of compact
support on $\partial X$, one has
\begin{equation}\label{5.11a}
\int_{\partial X}(|df|^{2} + f^{2}N(H)) \geq 0.
\end{equation}
Choose $f = f_{R,S}(r)$ such that $f = 1$ on $D(R) \subset \partial X =
{\mathbb R}^{2}$ and, for $r \geq R$, $f = (\log r - \log S)/(\log R - \log S)$,
for $S >> R >> 1$. One may choose $R$ and $S$ sufficiently large such that
$\int_{\partial X}|df|^{2} < \varepsilon$, for any given $\varepsilon > 0$.
This together with \eqref{5.10} implies \eqref{5.11}.

   It follows that $A = 0$ and hence by the Liouville theorem on ${\mathbb R}^{2}$,
$u = const$ on $\partial X$. Using the divergence constraint \eqref{5.4}, we also now
have $0 = \delta (A - H\gamma) = -u^{-1}D^{2}u(N, \cdot)$, and so $0 = D^{2}u(N, \cdot)
= dN(u) - A(du) = dN(u)$, so that $N(u) = const$.

  Thus, the full Cauchy data $(\gamma, u, A, N(u))$ for the static vacuum equations
is fixed and trivial: $\gamma$ is the flat metric, $A = 0$ and $u$, $N(u)$ are constant.
Observe that this data is realized by the family of flat metrics on $({\mathbb R}^{3})^{+}$
with either $u = const$ or $u$ equal to an affine function on $({\mathbb R}^{3})^{+}$.

  Suppose first
\begin{equation}\label{5.11ab}
u = const > 0 \ \ {\rm on} \ \  \partial X.
\end{equation}
The static vacuum equations \eqref{1.1} are then non-degenerate up to $\partial X$. The
unique continuation property for Einstein metrics with boundary, cf.~\cite{AH}, implies
that the Cauchy data uniquely determine the solution locally. Alternately, since the static
vacuum equations are non-degenerate up to $\partial X$ and since the boundary data
$(\partial X, \gamma, H)$ are real-analytic, elliptic regularity implies that the solution
$(M, g, u)$ is real-analytic up to $\partial M$. Such solutions are uniquely determined
(locally) by their Cauchy data. Hence, the limit $(X, g, u)$ is flat in this case.

  Moreover, the convergence to the limit is smooth everywhere. This again follows from
non-degeneracy and ellipticity. Briefly, the potential equation $\Delta u = 0$
gives a boost on the regularity of $u$, (given background regularity on $g$). One
substitutes this into the main static vacuum equation $uRic = D^{2}u$, giving thus
a boost to the regularity of $Ric$, inducing then a boost to the regularity of $g$.
This in turn further boosts the regularity of $u$ via the potential equation.
Bootstrapping gives $C^{m,\alpha}$ convergence, up to the boundary, in regions where
$u > 0$, (given that $u$ is $C^{m,\alpha}$ at the boundary). In sum, one has a
contradiction to \eqref{5.8}.

\smallskip

  Thus, suppose instead
\begin{equation}\label{5.11b}
u = 0 \ \ {\rm on} \ \ \partial X.
\end{equation}
This situation is more complicated. It is also more difficult to prove smooth convergence
in this situation (one may have $x$ in \eqref{5.8} at $\partial X$). Moreover, there are
in fact non-flat static vacuum solutions with flat Cauchy data as above with $u = 0$ on
$\partial X$, (so-called toroidal black holes, cf.~\cite{P} and \cite{Th}). Thus the
unique continuation results used above are false in this degenerate situation where the
boundary becomes characteristic.

   We observe first that the solution $(X, g)$ is still real-analytic up to $\partial X$
in this case. This follows since the 4-metric $g^{4} = u^{2}d\theta^{2} + g_{X}$ is Einstein
($Ric = 0$) and is $C^{1,\alpha}$ up to the horizon or vanishing locus $\partial X = \{u = 0\}$.
Elliptic regularity for the Einstein equations then implies that $g^{4}$ is real-analytic,
and hence so are $u$, $g_{X}$ up to $\partial X$.

  To see this in more detail, let $U$ be a chart neighborhood of $\partial X$ in $X$, so that
$U$ is diffeomorphic to a half-ball in ${\mathbb R}^{3}$ with boundary a disc $D^{2} \subset
{\mathbb R}^{2}$. Over each $p \in U$, one has a circle of length $2\pi u(p)$, with
$u \rightarrow 0$ as $p \rightarrow \partial X$. From the work above, $u$ is
$C^{1,\alpha}$ up to $\partial X$ with $N(u) = const$ at $\partial X$. Note that
$N(u) \neq 0$ at $\partial X$. For if $N(u) = 0$ at $\partial X$, since also $u = 0$
at $\partial X$ and $u$ is harmonic ($\Delta u = 0$), the unique continuation property for
harmonic functions implies that $u = 0$ in $X$, giving a contradiction. By rescaling $u$
if necessary, one may thus assume that $N(u) = 1$ at $\partial X$. This implies that the
4-metric $g^{4}$, defined on $B^{4}\setminus D^{2}$ extends to a $C^{1,\alpha}$ metric
on the 4-ball $B^{4}$. (The coordinate $\theta$ is an angular variable in ${\mathbb R}^{2}$
in polar coordinates, shrinking down to the origin on approach to $\partial X$). It is
well-known, cf.~\cite{Be}, that any $C^{1,\alpha}$ weak solution to the Einstein equations
is real-analytic (in harmonic or geodesic normal coordinates), which gives the claim
above.

  To prove the limit is in fact flat, and that one has strong convergence, we need to
use the outer-minimizing property again. Thus, first note that \eqref{5.9} holds everywhere
on the limit $(X, g)$ near $\partial X$, not just at $\partial X$; here $A$ is the
$2^{\rm nd}$ fundamental form of the level sets $S(t)$ of $t = dist(\partial X, \cdot)$,
etc. We have already established $N(H) = 0$ at $\partial X$, via the outer-minimizing
property and the corresponding stability of the $2^{\rm nd}$ variation operator
\eqref{5.11a}. Taking then the derivative of \eqref{5.9} in the normal direction
gives,
\begin{equation}\label{5.11c}
NN(H) = -\langle A', A \rangle - \langle A^{2}, A \rangle - H'H + {\tfrac{1}{2}}s_{\gamma}',
\end{equation}
where $A' = \nabla_{N}A$. At $\partial X$, the first three terms vanish while
$(s_{\gamma}')_{k} = -\Delta (tr k) + \delta \delta k - \langle Ric_{\gamma},
k \rangle = 0$, since $k = 2A = 0$. Thus
\begin{equation}\label{5.11d}
NN(H) = 0,
\end{equation}
at $\partial X$ and it follows that the $3^{\rm rd}$ variation of the area of
$\partial X$ in the unit normal direction vanishes.

  Now choose $f = f_{R,S}$ as following \eqref{5.11a} with $R$, $S$ large. Let
$S_{tf} = exp_{\partial X}(tf(x))$, where $exp_{\partial X}$ is the normal exponential
map of $\partial X$ into $X$. Letting $v(t) = area S_{tf}$, one has
\begin{equation}\label{5.12a}
v(t) = v(0) + {\tfrac{1}{2}}v''(0)t^{2} + {\tfrac{1}{6}}v'''(0)t^{3} +
{\tfrac{1}{24}}v^{''''}(0)t^{4} + O(t^{5}).
\end{equation}
The expansion \eqref{5.12a} is valid for all $t$ sufficiently small, $|t| \leq
\delta_{0}$, with $\delta_{0}$ independent of $R$, $S$, since the area and its
derivatives are integrals of local expressions, and the local geometry of $X$ is
uniformly bounded in a tubular neighborhood of radius 1 about $\partial X$.
By the $2^{\rm nd}$ variational formula \eqref{5.11a} and \eqref{5.11}, for any
given $\varepsilon > 0$, one has
$$v''(0) \leq \varepsilon,$$
for $R$, $S$ sufficiently large. For the same reasons via \eqref{5.11d},
$$v'''(0) \leq \varepsilon.$$
It follows then from the outer-minimizing property \eqref{5.7} and \eqref{5.12a} that
for $R$, $S$ sufficiently large, one must have
\begin{equation}\label{5.12b}
v^{''''}(0) \geq -\varepsilon,
\end{equation}
again for any $\varepsilon = \varepsilon(R, S) > 0$. Using the vanishing of the
lower order terms, one computes that \eqref{5.12b} gives
\begin{equation}\label{5.12}
\int_{\partial X}f^{4}NNN(H) - 6f^{2}\langle df \cdot df, A' \rangle \geq
-\varepsilon.
\end{equation}

   On the other hand, taking the normal derivative of \eqref{5.11c} gives
\begin{equation}\label{5.13}
NNN(H) = -\langle A', A' \rangle - (H')^{2} + {\tfrac{1}{2}}s_{\gamma}''.
\end{equation}
We have $H' = N(H) = 0$ and $(s_{\gamma}')_{2A} = 2\Delta H + 2\delta \delta A -
\langle Ric_{\gamma}, 2A \rangle$. For $s_{\gamma}''$, one has $(\Delta H)' =
\Delta' H + \Delta H' = 0$ and $\langle Ric_{\gamma}, A \rangle' =
\langle (Ric_{\gamma})', A \rangle + \langle (Ric_{\gamma}), A' \rangle = 0$.
So at $\partial X$, $\frac{1}{2}s_{\gamma}'' = \delta \delta A'$. It follows then
from \eqref{5.12}-\eqref{5.13} that for $f = f_{R,S}$ as above
\begin{equation}\label{5.14}
\int_{\partial X}-f^{4}|A'|^{2} + f^{4}\delta \delta A' -
6f^{2}\langle df \cdot df, A' \rangle \geq -\varepsilon.
\end{equation}
Integrating the second term by parts gives $\int \langle D^{2}f^{4}, A' \rangle =
\int \langle 4f^{3}D^{2}f + 12f^{2}\langle df\cdot df, A' \rangle$. Using the
Cauchy-Schwarz and Young inequalities,  \eqref{5.14} then implies, for any $\mu$
small,
$$\int_{\partial X}f^{4}|A'|^{2} \leq \mu \int_{\partial X}f^{4}|A'|^{2} +
 C\mu^{-1}\int_{\partial X}f^{2}|D^{2}f|^{2} + C\mu^{-1}\int_{\partial X}
|df|^{4} + \varepsilon.$$
Choosing $\mu$ small, the first term on the right may be absorbed into the left,
while simple computation shows that the last two terms become arbitrarily small
for $R$ and $S$ sufficiently large. It follows that
\begin{equation}\label{5.14a}
A' = 0
\end{equation}
and so $NNN(H) = 0$ at $\partial X$. The Riccati equation
\begin{equation}\label{5.14b}
A' = \nabla_{N}A = - A^{2} - R_{N},
\end{equation}
where $R_{N}(V,W) = \langle R_{g}(V,N)N, W \rangle$, thus gives $R_{N} = 0$ at
$\partial X$, and so via the Gauss and Gauss-Codazzi equations $R_{g} = 0$ at
$\partial X$. Thus the full ambient curvature vanishes at $\partial X$.

  One can now continue inductively in the same way to see that $A$ and $R$ vanish to
infinite order at $\partial X$. A simpler method proceeds as follows. The Riccati
equation \eqref{5.14b} holds along the level sets $S(t)$ of $t = dist(\partial X, \cdot)$.
Since $s_{g} = 0$ and $dim X = 3$, $R_{N} = -*Ric$, i.e.~$R_{N}(v,v) = -Ric(w,w)$, where
$(N, v, w)$ are an orthonormal basis. Via the static vacuum equations, this gives
$\nabla_{N}A = -A^{2} + u^{-1}*(D^{2}u)$. Rescale $u$ if necessary so that
$N(u) = 1$ at $\partial X$ and set $v = u - t$. Since $A = D^{2}t$,
one then obtains
\begin{equation}\label{5.15}
\nabla_{N}A = -A^{2} + u^{-1}(*A) + u^{-1}(*D^{2}v).
\end{equation}
This is a system of ODE's for $A$, singular at $\partial X = {\mathbb R}^{2}$, but
with indicial root 1. From the work above, we have $v = O(t^{2})$ and $A = O(t^{2})$.
Writing $A = t^{2}B$ and substituting in \eqref{5.15} shows that $A = O(t^{3})$.
Also, by the computation following \eqref{5.13}, $s_{\gamma} = O(t^{3})$ on $S(t)$,
and hence using the scalar constraint \eqref{5.3} and the relation $R_{NN} = u^{-1}NN(u)$,
this in turn implies $v = O(t^{3})$, and so on. It follows that $(g, u)$ agree with
a flat solution to infinite order at $\partial X$. Since the solution $(X, g, u)$ is
analytic up to $\partial X$, it follows that $(X, g, u)$ is flat, as claimed.

\medskip

  Next we claim that one has strong convergence to the limit, so that \eqref{5.5a} is
preserved in the limit, i.e.~\eqref{5.8} holds, contradicting the fact that the limit
is flat. Note first that if $x$ in \eqref{5.8} is in the interior of $X$, then
strong ($C^{\infty}$) convergence is immediate, by the interior estimates
\eqref{5.2}, i.e.~their higher derivative analogs. Thus we may assume that
$x \in \partial X$.

   From the work above, we know that $|R|$ is uniformly bounded everywhere on
$(M, g_{i}')$ and $|R| \rightarrow 0$ everywhere away from $\partial M \rightarrow
\partial X$, so $|R|$ jumps quickly from 1 to 0 near $x_{i}$. The main point is
to prove that
\begin{equation}\label{5.16}
|R|(x) \rightarrow 0,
\end{equation}
for all $x \in \partial M \rightarrow \partial X$; it is then easy to prove that
$|R| \rightarrow 0$ on $(M, g_{i}')$, cf.~\eqref{5.18} below. To prove \eqref{5.16},
note that the estimates \eqref{5.9}-\eqref{5.14a} above at $\partial X$ also hold
on the blow-up sequence at $\partial M$, (since $H \rightarrow 0$ and $s_{\gamma}
\rightarrow 0$ on $\partial M$). It follows then by these arguments that for any
$R < \infty$ and $D(R)$ the $R$-ball about any base point $y_{i} \in (\partial M,
\gamma_{i}')$ converging to $y \in \partial X$,
\begin{equation}\label{5.16a}
\int_{D(R)}|A|^{2} + |\nabla_{N}A|^{2} \rightarrow 0.
\end{equation}
To proceed further, consider the divergence constraint $\delta A = -dH - Ric(N, \cdot)$
on $(\partial M, \gamma_{i}')$, as in \eqref{5.4}. The equations $\delta A =
\chi_{1}$, $d tr A = \chi_{2}$ form an elliptic first order system in 2-dimensions,
and so one has elliptic estimates. Since $H \rightarrow 0 \in C^{m-1,\alpha}$ on
$(\partial M, \gamma_{i}')$, $dH \rightarrow 0$ in $C^{m-2,\alpha}$. Also, by
assumption \eqref{5.5a}, $Ric(N, \cdot)$ is bounded in $L^{\infty}$. It follows
then from elliptic regularity that $A$ is bounded in $L^{1,p}$. By the scalar
constraint \eqref{5.3}, $R_{NN}$ is then also bounded in $L^{1,p}$ and since
$tr R_{N} = R_{NN}$, it follows from \eqref{5.16a} and \eqref{5.14b} that
\begin{equation}\label{e5.17}
A \rightarrow 0 \ \ {\rm and} \ \ R_{NN} \rightarrow 0 \ \ {\rm in} \ \
C_{loc}^{\alpha}(\partial M).
\end{equation}

  Next the Einstein equation $Ric_{g_{\mathcal{M}}} = 0$ on $(\mathcal{M}, g_{\mathcal{M}})$ implies that
$\delta_{\mathcal{M}} R_{g_{\mathcal{M}}} = 0$. Hence
\begin{equation}\label{5.171}
0 =  (\delta_{g_{\mathcal{M}}}R)(N, N, \cdot) = \delta_{g_{\mathcal{M}}}(R(\cdot, N)N, \cdot)
+ 2R(e_{\alpha}, \nabla_{e_{\alpha}}N)N = \delta_{g_{\mathcal{M}}}(R(\cdot, N)N, \cdot),
\end{equation}
since one may choose a basis in which $\nabla_{e_{\alpha}}N = A(e_{\alpha}) =
\lambda_{\alpha}e_{\alpha}$. Let $V$ be the unit vertical vector, and note that
$\nabla_{V}V = -d\nu$, where $\nu = \log u$. Then $\nabla_{V}(R(V, N)N, \cdot) =
R(d\nu, N)N$, while $\nabla_{N}(R(N, N)N, \cdot) = 0$. Hence, for $R_{N}$ as
in \eqref{5.14b}, these computations on $\partial M$ give
$$\delta R_{N} = -R_{N}(d\nu),$$
where the divergence $\delta$ and $R_{N}$ are taken on $(\partial M, \gamma_{i}')$.
Since $R$ is bounded in $L^{\infty}$ and $tr R_{N} = R_{NN}$ is bounded in
$L^{1,p}$, elliptic regularity gives
\begin{equation}\label{5.17a}
||R_{N}||_{L^{1,p}} \leq C||d\nu||_{L^{p}},
\end{equation}
again on compact domains in $\partial M$ converging to a compact domain in
$\partial X$.

  To control $d\nu$ in \eqref{5.17a}, recall that by \eqref{5.5a} the ambient curvature
$R$ is bounded, and so $R$ restricted to $\partial M$ is also bounded. Via the
static vacuum equations \eqref{1.1}, this implies that
\begin{equation}\label{5.17b}
N(\nu)A + u^{-1}D^{2}u
\end{equation}
is bounded on $(\partial M, \gamma_{i}')$, where $D^{2}u$ is the Hessian of
$u|_{\partial M}: \partial M \rightarrow {\mathbb R}^{+}$. Now we claim that
each term in \eqref{5.17b} is bounded, i.e.~there exists $K$ such that
\begin{equation}\label{5.17c}
|u^{-1}D^{2}u| \leq K, \ \ |N(\nu)A| \leq K,
\end{equation}
pointwise, on domains converging to a bounded domain in $\partial X$. To prove
\eqref{5.17c}, suppose instead that $|u^{-1}D^{2}u| \rightarrow \infty$ at some
sequence of base points $y_{i} \rightarrow y \in \partial X$. Without loss of
generality, we may assume the points $y_{i}$ realize the maximum of $|u^{-1}D^{2}u|$
on $D_{y_{i}}(10) \subset (\partial M, \gamma_{i}')$ (possibly up to a factor of
2), where $D_{y_{i}}(10)$ is the geodesic disc of radius $10$ in $\partial M$ centered at $y_{i}$ . As before, one may then rescale the metrics $g_{i}'$ further to $g_{i}''$ so that
$|u^{-1}D^{2}u|(y_{i}) = 1$ and hence the full curvature $R \rightarrow 0$ in this
scale. Also, renormalize $u$ if necessary so that $u(y_{i}) = 1$. Note that
$|du|(y_{i})$ must be bounded. For if $|du|(y_{i})$ were too large, it
follows, (e.g.~by a still further rescaling), that $u$ would be close to an
affine function on ${\mathbb R}^{2}$ and hence $u$ would assume negative values
in bounded distance to $y_{i}$. Since $u > 0$ everywhere, this is impossible.
Thus, by integration along paths, $u$ is bounded in $L^{1,\infty}$ in the scale
$g_{i}''$, within bounded distance to $y_{i}$.

  Now working in the scale and normalization above, from divergence constraint
\eqref{5.4}, one has $-u\delta(A - H\gamma) = dN(u) - A(du)$. Since $A$ is bounded
in $L^{1,p}$ and $u$ and $du$ are bounded in $L^{\infty}$, it follows that $dN(u)$
is bounded in $L^{p}$. Thus, $N(u) = c + \phi$, where $\phi$ is bounded in $L^{1,p}$.
Here $c$ is a constant which may, and in fact does, go to $\pm \infty$, in the
$u$-normalization $u(y_{i}) = 1$ above.

  The trace equation \eqref{5.5} in this scale and normalization gives
\begin{equation}\label{5.17d}
\Delta u + H(c + \phi) = uf,
\end{equation}
where $f$ is bounded in $L^{1,p}$. Also, $Hc$ is bounded, since $N(u)A$ is bounded
via \eqref{5.17b}, and so $Hc \rightarrow c'$, for some constant $c'$ on $\partial X$.
Since $\phi$ is also bounded in $L^{1,p}$, it follows that $u$ is bounded in $L^{3,p}$,
and hence (in a subsequence), $u$ converges in $C^{2,\alpha}$ to its limit on
$\partial X = {\mathbb R}^{2}$. Hence $D^{2}u$ converges in $C^{\alpha}$ to its limit
on ${\mathbb R}^{2}$. On the limit, since $H \rightarrow 0$, the trace equation
\eqref{5.17d} becomes
\begin{equation}\label{5.17e}
\Delta u + c' = 0.
\end{equation}
If $c' = 0$ then since $u > 0$, $u = const$ and hence $D^{2}u = 0$, giving a
contradiction. If $c' \neq 0$, then again since $u > 0$, one must have $c' < 0$
and $u$ is a quadratic polynomial on ${\mathbb R}^{2}$. Since $u$ is harmonic
on the blow-up sequence and $c' < 0$ implies $N(u) \rightarrow -\infty$, this is
also easily seen to be inconsistent with the requirement $u > 0$ everywhere.
This proves \eqref{5.17c} holds.

  Returning to \eqref{5.17a}, we may work in the normalization above that
$u(y_{i}) = 1$ and then \eqref{5.17c} implies that $d\nu = u^{-1}du$ is bounded
in $L^{\infty}$ in bounded domains about $y_{i}$. Hence by \eqref{5.17a}, $R_{N}$
is bounded in $L^{1,p}$ and so bounded in $C^{\alpha}$. Since $R_{N} \rightarrow 0$
in $L^{2}$ locally on $\partial M$, one has
\begin{equation}\label{5.17f}
R_{N} \rightarrow 0 \ \ {\rm in} \ \ C_{loc}^{\alpha}(\partial M).
\end{equation}
This is the main part of the estimate \eqref{5.16}.

  Next, one needs the same result for the $Ric(N, T)$ term. To do this, take the normal
derivative of the scalar constraint \eqref{5.5}, to obtain
\begin{equation}\label{5.17g}
\Delta N(u) + \Delta'u + N(H)N(u) + HNN(u) = {\tfrac{1}{2}}N(u)[|A|^{2} - H^{2} +
s_{\gamma}] + {\tfrac{1}{2}}uN[|A|^{2} - H^{2} + s_{\gamma}].
\end{equation}
Since $\gamma' = 2A$, a standard formula for the variation of the Laplacian,
(cf.~\cite{Be}), gives $\frac{1}{2}\Delta'u = -\langle D^{2}u, A \rangle +
\langle du, \beta (A) \rangle$ which is bounded in $L^{\infty}$. Also the terms $N(H)$,
$H$ and $NN(u) = -\frac{u}{2}[|A|^{2} - H^{2} + s_{\gamma}]$ all go to 0 in $L^{\infty}$.
For the right side of \eqref{5.17g}, the coefficient of the first term goes to 0 in $L^{\infty}$
while via \eqref{5.17f} above, $uN[|A|^{2} - H^{2} + s_{\gamma}]$ is bounded in
$L^{p}$ provided $N(s_{\gamma})$ is, and this in turn follows from an $L^{p}$ bound on
$\delta \delta A$.

     To obtain this, return to \eqref{5.171} but with $(N, X)$ in place of
$(N,N)$, with $X$ tangent to the boundary. Arguing as above, it follows that
$\delta(R(\cdot, N)X, \cdot)$ is bounded in $L^{p}$. Via the Gauss-Codazzi equations
$dA(X,Y,Z) = \langle R(N,X)Y, Z\rangle$, it follows that $\delta dA$ is bounded in $L^{p}$.
Similarly, in the normalization above, $\delta A = D^{2}u(N, \cdot) = \nabla_{N}du$ modulo
lower order terms. Taking the exterior derivative $d$ of this, one easily obtains that $d\delta A$
is bounded in $L^{p}$. This shows that $\Delta A$ is bounded in $L^{p}$, which of course
gives the same for $\delta \delta A$ by elliptic regularity.

   Finally, as in the proof of \eqref{5.17c}, $N(u) = c + \phi$ with $\phi$ bounded
in $L^{1,p}$. Thus it follows from the divergence constraint as preceding \eqref{5.17d}
and elliptic regularity that $dN(u)$ is bounded in $L^{1,p}$ and so converges in
$C^{\alpha}$ on $\partial M$. Since
$$-uRic(N, \cdot) = dN(u) - A(du),$$
it now follows that $Ric(N, \cdot)$ converges to its limit, necessarily 0, in
$C^{\alpha}$ on $\partial M$. Lastly, $\delta A = -dH - Ric(N, \cdot) \rightarrow 0$
in $C^{\alpha}$ and hence by elliptic regularity, $A \rightarrow 0$ in $C^{1,\alpha}$
so that again via the Gauss-Codazzi equations $dA(X,Y,Z) = \langle R(N,X)Y, Z\rangle
\rightarrow 0$ in $C^{\alpha}$. Combining the computations above now proves the
estimate \eqref{5.16}.

  To complete the proof, we have $|R| \rightarrow 0$ on $\partial M \rightarrow
\partial X$. On the 4-manifold $\mathcal{M} = M\times_{u}S^{1}$, the Einstein equations give
the
inequality
\begin{equation}\label{5.18}
\Delta_{\mathcal{M}} |R_{g_{\mathcal{M}}}| + c|R_{g_{\mathcal{M}}}|^{2} \geq 0,
\end{equation}
where $c$ is a fixed numerical constant, $\Delta_{\mathcal{M}}$ is the 4-Laplacian and
$R_{g_{\mathcal{M}}}$ is the curvature tensor on $\mathcal{M}$. One has $|R_{g_{\mathcal{M}}}| = |R|$, (up to a
constant). We have proved above that $R \rightarrow 0$ in $L^{2}$ locally on
$(\mathcal{M}, (g_{\mathcal{M}})_{i})$ and $R \rightarrow 0$ pointwise at $\partial \mathcal{M}$. It follows
from the deGiorgi-Nash-Moser estimates for domains with boundary,
(cf.~[GT, Thm.~8.25]), that $R \rightarrow 0$ pointwise on $\mathcal{M}$ and hence on
$M$, contradicting \eqref{5.5a}. This completes the proof.

{\endproof}

  Next, we show that the potential function $u$ is also controlled by the
boundary data of $\Pi_{B}$.

\begin{corollary}\label{c5.2}
For $(M, g, u) \in {\mathcal E}^{o}$, there is a constant $U_{0}$ depending
only on the boundary data $(\gamma, H) \in \Pi^{o}({\mathcal E}^{o})$ such that
\begin{equation}\label{5.20}
u \leq U_{0},
\end{equation}
on $M$. Moreover, if $H \geq H_{0} > 0$ on $\partial M$, then there exists $U_{0}$
as above depending in addition only on $H_{0}$, such that on $M$,
\begin{equation}\label{5.21}
u \geq U_{0}^{-1}.
\end{equation}
\end{corollary}

{\bf Proof:} Let $S(s) =  \{x \in (M, g): dist(x, \partial M) = s\}$ be the geodesic
'sphere' about $\partial M$. Choose a fixed base point $x_{0} \in S(1)$ and
suppose one has the bound
\begin{equation}\label{5.22}
c_{0}^{-1} \leq u(x_{0}) \leq c_{0}.
\end{equation}
By Theorem 4.1, the geometry of the annular region $A(\frac{1}{2}, 2)$ about $S(1)$
is uniformly controlled by the boundary data $(\gamma, H)$ and so by integration
of the static vacuum equations $uRic = D^{2}u$ along paths in $A(\frac{1}{2},2)$,
one has
\begin{equation}\label{5.23}
C_{0}^{-1} \leq u \leq C_{0},
\end{equation}
in $A(\frac{3}{4}, \frac{3}{2})$, where $C_{0}$ depends only on $c_{0}$ and
$(\gamma, H)$.

  To remove the dependence of $C_{0}$ in \eqref{5.23} on $c_{0}$ in \eqref{5.22}, we
need better control on the large-scale behavior of $u$. To do this, it is proved in
[An2, Lemma 3.6], that there is a constant $K$, depending only on $C_{0}$, such
that for all $s \geq 1$,
\begin{equation}\label{5.25}
\sup_{S_{c}(s)}|du| \leq K(v_{c}(s))^{-1},
\end{equation}
where $v_{c}(s) = area S_{c}(s)$ and $S_{c}(s)$ is any component of the geodesic
sphere $S(s)$. We point out that \eqref{5.25} holds for general static vacuum
solutions, not only those in ${\mathcal E}^{o}$ for instance. The estimate \eqref{5.25}
is proved by studying the behavior of the harmonic potential $\log u$ on the Ricci-flat
4-manifold $(\mathcal{M}, g_{\mathcal{M}}, \log u)$ and then reducing to $(M, g_{M}, u)$.

  Consider now the conformally equivalent metric
\begin{equation}\label{5.25a}
\widetilde g = u^{2}g.
\end{equation}
It is well-known that the static vacuum Einstein equations \eqref{1.1} are equivalent
to the equations $\widetilde Ric = 2(d\nu)^{2} \geq 0$, $\Delta_{\widetilde g}
\nu = 0$, $\nu = \log u$. The metric $\widetilde g$ thus has non-negative Ricci
curvature with harmonic potential $\nu$. These are exactly the properties used to
prove \eqref{5.25}, and a brief examination of its proof shows that \eqref{5.25}
also holds with respect to $\widetilde g$, i.e.
\begin{equation}\label{5.26}
\sup_{\widetilde S_{c}(s)}|d\nu|_{\widetilde g} \leq K(\widetilde v_{c}(s))^{-1},
\end{equation}
again with $K = K(C_{0})$. Since $(M, g)$ is asymptotically flat and $u \rightarrow const$
at infinity, the area growth of geodesic spheres $\widetilde v(s)$ in $(M, \widetilde g)$
satisfies $\widetilde v(s) / s^{2} \rightarrow \omega_{2}$, where $\omega_{2} = area S^{2}(1)$.
It follows then from the volume comparison theorem for Ricci curvature, (cf.~\cite{Pe}),
that
\begin{equation}\label{5.27}
\widetilde v(s) \geq \omega_{2}s^{2},
\end{equation}
for all $s \geq 1$. As above, by integration of \eqref{5.26} along a geodesic ray
starting from a suitable base point $x_{1} \in S(1)$ out to infinity, one sees that
\eqref{5.20} holds then globally on $M\setminus B(1)$, with $U_{0}$ again depending
only on $c_{0}$ in \eqref{5.22}. Using the static vacuum equations, the same
integration along paths gives such a bound within $B(1)$. Thus, we see that
\eqref{5.20} follows from \eqref{5.22}.

  To prove \eqref{5.22}, suppose one has a static vacuum solution $(M, g, u)$
with
\begin{equation}\label{5.28}
u(x_{0}) = \varepsilon.
\end{equation}
Renormalize $u$ to $\bar u = u/u(x_{0})$, so that $\bar u(x_{0}) = 1$ and
$\bar u \rightarrow \varepsilon^{-1}$ at infinity. Then \eqref{5.22} holds, and
hence so does \eqref{5.26}-\eqref{5.27}. Again by integration along geodesics
starting at $x_{0}$ and diverging to infinity, it follows that
$$u \leq U_{1},$$
where $U_{1}$ depends only on the boundary data of $\Pi_{B}$. This proves the lower
bound in \eqref{5.22}.

  To prove the upper bound, suppose instead
\begin{equation}\label{5.29}
u(x_{0}) = \varepsilon^{-1}.
\end{equation}
Then again we renormalize $u$ to $\bar u$ as above, so that now $u \rightarrow
\varepsilon$ at infinity. This does not directly give a lower bound on $\varepsilon$
via \eqref{5.26}-\eqref{5.27} as above. However, one may proceed as follows. First,
it is well-known that static vacuum solutions come in ``dual'' pairs, in that if
$(M, g, \bar u)$ is a static vacuum solution, then so is $(M, \hat g, \hat u)$
with $\hat g = \bar u^{4}g$, $\hat u = \bar u^{-1}$, cf.~\cite{A2} for instance. Then
\eqref{5.26}-\eqref{5.27} hold for $(M, \hat g, \hat u)$ which as before by
integration gives an upper bound $\hat u \leq U_{1}$ at infinity. Since near
infinity, $\hat u \simeq \varepsilon^{-1}$, this again gives a bound on
$\varepsilon^{-1}$. This completes the proof of \eqref{5.20}.

  To prove \eqref{5.21}, note that \eqref{5.22} has been proved above, and hence by
the maximum principle and normalization $u \rightarrow 1$ at infinity, \eqref{5.21}
holds in the exterior region $M \setminus B(1)$. Thus, one only needs to consider
the behavior near $\partial M$. For this, suppose $(M, g, u)$ is a static vacuum
solution, $C^{2,\alpha}$ up to $\partial M$ with $u \geq 0$ on $\bar M =
M \cup \partial M$. If $H \geq H_{0} > 0$ on $\partial M$, we claim that necessarily
$$u > 0 \ \ {\rm on} \ \ \partial M.$$
For if $u = 0$ at some point $z\in \partial M$, then by \eqref{5.5},
$\Delta_{\partial M}u + HN(u) = 0$ at $z$. Since $0 = u(z)$ is a global minimum
for $u$, one has $\Delta_{\partial M}u \geq 0$ and by the Hopf maxmimum principle,
$N(u) > 0$ at $z$. This gives a contradiction if $H > 0$. The same arguments
prove the existence of a lower bound \eqref{5.21} by a contradiction argument,
taking a sequence and passing to a limit, using Theorem 4.1.

{\endproof}

  The previous results now lead quite easily to the following main result of this
section.

\begin{corollary}\label{c5.3}
The boundary map
$$\Pi^{o}: {\mathcal E}^{o} \rightarrow Met^{m,\alpha}(\partial M) \times C_{+}^{m-1,\alpha}
(\partial M),$$
is almost proper.
\end{corollary}

{\bf Proof:} Let $(M, g_{i}, u_{i})$ be a sequence of static vacuum solutions in
${\mathcal E}^{o}$, with $\Pi^{o}(g_{i}, u_{i}) = (\gamma_{i}, H_{i})$.
Supposing $(\gamma_{i}, H_{i}) \rightarrow (\gamma, H)$ in
$Met^{m,\alpha}(\partial M) \times C_{+}^{m-1,\alpha}(\partial M)$, we need to
prove that the sequence $(g_{i}, u_{i})$ has a subsequence converging in
$C^{m,\alpha}(M)$, modulo diffeomorphisms, to a limit $(M, g, u) \in {\mathcal E}$.

   The curvature bound \eqref{5.1} and control of the intrinsic and extrinsic
geometries of the boundary metrics first implies the metrics $g_{i}$ cannot
collapse within bounded distance to $\partial M$, i.e.~there is a fixed constant
$i_{0} > 0$ such that the injectivity radius of $(M, g_{i})$ satisfies
$$inj_{g_{i}}(x) \geq i_{0},$$
for $dist_{g_{i}}(x, \partial M) \leq K$. By the compactness theorem in, for
instance \cite{AT}, it follows that a subsequence of $(M, g_{i})$ converges in
$C^{m,\alpha}$, (and $C^{\infty}$ in the interior), uniformly on bounded domains
containing $\partial M$, to a limit $(M', g)$. One has $\partial M' = \partial M$
and $g$ is a complete Riemannian metric on $M'$, $C^{m,\alpha}$ up to $\partial M$
and $C^{\infty}$ in the interior.

  By Corollary 4.2, the potential functions $u_{i}$ also converge in $C^{m,\alpha}$,
(in a subsequence) to a limit potential function $u$ on $M'$, and the pair $(g, u)$
gives a solution of the static vacuum Einstein equations. Since $u = \lim u_{i}$,
it follows that
\begin{equation}\label{5.30}
u \leq U_{0},
\end{equation}
on $M'$, for $U_{0}$ as in \eqref{5.20}. Clearly the boundary metric and mean curvature
of $(M', g, u)$ are given by the limit values $(\gamma, H)$. To prove that $(M', g, u)
\in {\mathcal E}$, one then needs to prove that $(M', g, u)$ is asymptotically
flat and $M'$ is diffeomorphic to $M$. Note that since the convergence above is only
uniform on compact sets, apriori there need not be any relation between the asymptotic
structure of $(M', g, u)$ and $(M, g_{i}, u_{i})$ for any given $i$.

  The equation \eqref{5.25} holds on each $(M, g_{i}, u_{i})$ and by Corollary
4.2, the constant $K$ is uniform, independent of $i$. Moreover, as in the
proof of Corollary 4.2, there is a geodesic ray $\sigma = \sigma_{i}$ starting
at any fixed base point in $S(1)$ and diverging to infinity, such that on the
component $S_{c}(s)$ of $S(s)$ containing $\sigma$, one has
$$\sup_{S_{c}(s)}|du_{i}| \leq Ks^{-2},$$
where $K$ is independent of $i$. Since $u_{i}$ is harmonic, by elliptic regularity,
(and scaling), a similar estimate holds for higher derivatives of $u_{i}$, and via
the static vacuum equations, it follows that
$$\sup_{S_{c}(s)}|R_{g_{i}}| \leq Cs^{-3},$$
with again $C$ independent of $i$. This means that the metrics $(M, g_{i}, u_{i})$
become asymptotically flat at infinity uniformly, at a rate independent of $i$.
For $R$ sufficiently large, $R \geq R_{0}$ independent of $i$, the geodesic spheres
$S(R)$ and annuli $A(R, 2R)$ are close to Euclidean spheres and annuli, (when scaled
by $R^{-1}$), and hence the geometry is close to that of Euclidean space; there can
be no branching or joining of different components of $S(R)$ for $R \geq R_{0}$. This
implies that the limit $(M', g, u)$ has a single asymptotically flat end, and $M'$
is diffeomorphic to $M$.

{\endproof}

\begin{remark}\label{r5.4}
{\rm The results above also show that the boundary map $\Pi_{B}$ is almost proper not
only on ${\mathcal E}^{o}$ but also its closure $\overline{\mathcal E}^{o}$.
In other words, if $(M, g_{i}, u_{i})$ is a sequence of static vacuum solutions in
${\mathcal E}^{o}$ with boundary data $(\gamma_{i}, H_{i})$ and $(\gamma_{i}, H_{i})
\rightarrow (\gamma, H)$ in $Met^{m,\alpha}(\partial M)\times C^{m-1,\alpha}(\partial M)$,
then $(M, g_{i}, u_{i})$ converges in $C^{m,\alpha}$ (in a subsequence) to a limit $(M, g, u)$ in
$\overline{\mathcal E}^{o}$, with $H \geq 0$.

  To see this, note that the constant $\Lambda$ in Theorem 4.1 does not depend on a
positive lower bound on $H$, and so \eqref{5.1}-\eqref{5.1a} hold for the sequence
$(M, g_{i}, u_{i})$ above. Of course we are using here the fact that $\partial M$ is
outer-minimizing on the sequence $(M, g_{i}, u_{i})$. Similarly in Corollary 4.2,
the upper bound $U_{0}$ on $u_{i}$ does not depend on a lower bound for $H$. One may
then use the argument concerning \eqref{5.28} to show that $u_{i}$ cannot go to 0
on $M$ away from $\partial M$. The proof of Corollary 4.3 also does not require a
bound on $H$ away from 0.

  It is also worth pointing out a brief examination of the proof shows that the results of
this section only require that $\partial M$ is outer-minimizing in a neighborhood of
arbitrarily small but fixed size (depending on $(M, g, u)$) about $\partial M$. }
\end{remark}

\section{Degree of $\Pi^{P}$.}
\setcounter{equation}{0}

  By \cite{Sm}, a smooth proper Fredholm map $F: B_{1} \rightarrow B_{2}$ of Fredholm
index 0 between connected Banach manifolds $B_{1}$, $B_{2}$ has a well-defined degree
(mod 2). Namely, if $y$ is a regular value of $F$, then $F^{-1}(y)$ is a finite set of
points, and $deg_{{\mathbb Z}_{2}}F$ is just the cardinality of $F^{-1}(y)$ (mod 2).
In fact if $B_{1}$ and $B_{2}$ are oriented, then $F$ has a well-defined degree
in ${\mathbb Z}$, cf.~\cite{ET}, \cite{BFP}.

   From the discussion in \S 1, the boundary map $\Pi_{B}: {\mathcal E}^{+} \rightarrow
Met(\partial M)\times C_{+}(\partial M)$, although Fredholm, is not proper. However,
by Theorem 1.2, the restricted boundary map $\Pi^{o}: {\mathcal E}^{o} \rightarrow
Met(\partial M)\times C_{+}(\partial M)$ is almost proper. This in fact suffices to obtain
a ${\mathbb Z}_{2}$-valued degree on natural domains within ${\mathcal E}^{o}$.

   Thus, as in the Introduction, let $\partial {\mathcal E}^{o}$ be the
boundary of ${\mathcal E}^{o}$ within the space ${\mathcal E}_{S}$ of static vacuum solutions.
Hence, $(M, g, u) \in \partial {\mathcal E}^{o}$ if and only if $(M, g, u)$ satisfies \eqref{1.12}
but not \eqref{1.10}. Let $Z = \Pi_{B}(\partial {\mathcal E}^{o}) \subset Met^{m,\alpha}(\partial M)\times
C^{m-1,\alpha}(\partial M)$ be the image of $\partial {\mathcal E}^{o}$ under the boundary
map $\Pi_{B}$ and set
$${\mathcal E}^{P} = (\Pi^{o})^{-1}([Met^{m,\alpha}(\partial M)\times
C^{m-1,\alpha}(\partial M)]\setminus Z).$$
Then, by Theorem 1.2 and construction, the induced boundary map
\begin{equation}\label{6.1}
\Pi^{P}: {\mathcal E}^{P} \rightarrow  [Met^{m,\alpha}(\partial M)\times
C^{m-1,\alpha}(\partial M)]\setminus Z
\end{equation}
is proper. In particular ${\mathcal E}^{P}$ has only finitely many components ${\mathcal E}^{P_{i}}$
and on each component the ${\mathbb Z}_{2}$-valued Smale degree is well-defined.

\medskip

  This discussion leads to Theorem 1.3. Let ${\mathcal E}^{P_{0}}$ be the
component of ${\mathcal E}^{P}$ containing the standard round flat solution,
equal to the exterior of the round ball in ${\mathbb R}^{3}$, with
$$\Pi^{P_{0}}: {\mathcal E}^{P_{0}} \rightarrow {\mathcal T}_{0}.$$

\begin{theorem}\label{t6.2}
The degree of $\Pi^{P_{0}}$ satisfies
\begin{equation}\label{6.4}
deg_{{\mathbb Z}_{2}} \Pi^{P_{0}} = 1.
\end{equation}
\end{theorem}

{\bf Proof:} The proof is based on the black hole uniqueness theorem \cite{I}, \cite{R},
\cite{BM}, that the Schwarzschild metrics
\begin{equation}\label{6.5}
g_{Sch}(m) = (1 - \frac{2m}{r})^{-1}dr^{2} + r^{2}g_{S^{2}(1)}, \ \
u = \sqrt{1 - \frac{2m}{r}},
\end{equation}
$r \geq 2m$, are the unique AF static vacuum metrics with a smooth horizon ${\mathcal H} =
\{u = 0\}$. The induced metric on $\partial M$ is $S^{2}(2m)$ - the round metric $\gamma_{2m}$
of radius $2m$ on $S^{2}$. The mean curvature satisfies $H = 0$. Of course the Schwarzschild
metrics are not in ${\mathcal E}_{S}$, but instead lie at the boundary $\partial {\mathcal E}_{S}$.

  Consider any sequence $\{(g_{i}, u_{i})\} \in {\mathcal E}^{P_{0}}$ for which
$\Pi^{P_{0}}(g_{i}, u_{i}) = (\gamma_{i}, H_{i}) \rightarrow (\gamma_{2m}, 0)$ smoothly. Clearly,
$\{(g_{i}, u_{i})\}$ is a divergent sequence in ${\mathcal E}^{o}$. By Corollary 4.3 and
Remark 4.4, a subsequence of $\{(g_{i}, u_{i})\}$ converges smoothly to a static
vacuum limit $(M, g, u)$. (Of course one may have $u = 0$ on $\partial M$). On this limit,
$$H = 0,$$
at $\partial M$, so that $\partial M$ is a minimal surface. From \eqref{5.5} one has
$2\Delta_{\partial M}u = u(|A|^{2} + s_{\gamma}) \geq 0$, and hence it follows from
the maximum principle that $u = 0$ on $\partial M$. Via the static vacuum equations
\eqref{1.1}, this implies further that $A = 0$ and $N(u) = const$ at $\partial M$.
The black hole uniqueness theorem then implies that any such limit is the Schwarzschild
metric, and so unique up to scaling. Thus one has uniqueness for the boundary data
$(\gamma, 0)$, so that almost all boundary metrics $\gamma$ cannot be realized with
$H = 0$ at $\partial M$, (the no-hair result).

  Given this background, suppose
\begin{equation}\label{6.6}
deg_{{\mathbb Z}_{2}} \Pi^{P_{0}} = 0.
\end{equation}
Then for any regular value $(\gamma, H)$ of $\Pi^{P_{0}}$, the finite set
$(\Pi^{P_{0}})^{-1}(\gamma, H)$, if non-empty, consists of at least two distinct static
vacuum solutions $(g^{1}, u^{1})$, $(g^{2}, u^{2})$. The regular values of $\Pi^{P_{0}}$
are generic (of second category) in the range space, by the Sard-Smale theorem. Choose
then a sequence of regular values $(\gamma_{i}, H_{i}) \rightarrow (\gamma_{+1}, 0)$
smoothly. (We set $m = 1/2$ here). Suppose for now that $(\Pi^{P_{0}})^{-1}(\gamma_{i}, H_{i})$
is non-empty; this will proved to be the case later.

   Let $(g_{i}^{1}, u_{i}^{1})$, $(g_{i}^{2}, u_{i}^{2})$ be any pair of corresponding
distinct sequences in $(\Pi^{P_{0}})^{-1}(\gamma_{i}, H_{i})$. By Corollary 4.3 and Remark 4.4,
the sequences $(g_{i}^{1}, u_{i}^{1})$, $(g_{i}^{2}, u_{i}^{2})$ have $C^{m,\alpha}$
convergent subsequences to limits $(g_{\infty}^{1}, u_{\infty}^{1})$, $(g_{\infty}^{2},
u_{\infty}^{2})$ in $\overline{\mathcal E}^{P_{0}}$ and by the uniqueness above
$$g_{\infty}^{1} = g_{\infty}^{2} = g_{Sch}(m),$$
with $m = 1/2$, with $u_{\infty}^{1} = u_{\infty}^{2} = u$ in \eqref{6.5}.

   This implies that near $g_{Sch}$, the boundary map $\Pi^{P_{0}}$ is not locally 1-1,
and so presumably $D\Pi_{B}$ has a non-trivial kernel at $g_{Sch}$. (Note however
that $g_{Sch} \notin {\mathcal E}_{S}$). We claim this is impossible. To prove the claim,
let
$${\sf g}_{Sch} = u^{2}d\theta^{2} + g_{Sch},$$
be the 4-dimensional Schwarzschild metric on ${\mathbb R}^{2}\times S^{2}$, and similarly let
$${\sf g}_{i}^{j} = (u_{i}^{j})^{2}d\theta^{2} + g_{i}^{j},$$
be the 4-dimensional static Ricci-flat metrics associated to $(g_{i}^{j}, u_{i}^{j})$.
By Lemma 2.2, without loss of generality we may assume that each ${\sf g}_{i}^{j}$ is
in Bianchi gauge with respect to ${\sf g}_{Sch}$, so that, as in \eqref{2.10}-\eqref{2.11},
$$\beta_{{\sf g}_{Sch}}({\sf g}_{i}^{j}) = 0,$$
for $j = 1, 2$ and $i$ sufficiently large. By the smoothness of the convergence above,
one may write
\begin{equation}\label{6.7}
{\sf g}_{i}^{j} = {\sf g}_{Sch} + \varepsilon_{i}^{j}\kappa_{i}^{j} +
O((\varepsilon_{i}^{j})^{2}),
\end{equation}
where $L(\kappa_{i}^{j}) = 0$ and $L$ is the linearized Einstein operator \eqref{2.6} at
${\sf g}_{Sch}$. The data ${\sf g}_{i}^{j}$, ${\sf g}_{Sch}$ and $\kappa_{i}^{j}$
are all smooth, (up to the boundary). The forms $\kappa_{i}^{j}$ are only unique
up to multiplicative constants, which will be determined by choosing $\varepsilon_{i}^{j}$
so that the $C^{1,\alpha}$ norm of ${\sf g}_{i}^{j} - {\sf g}_{Sch}$ equals
$\varepsilon_{i}^{j}$. Thus the $C^{1,\alpha}$ norm of $\kappa_{i}^{j}$ is on the
order of 1. Note that $\kappa_{i}^{j}$ decays to 0 at infinity, so it is basically
supported within compact regions of $M$. Let $\varepsilon_{i} = \max(\varepsilon_{i}^{1},
\varepsilon_{i}^{2})$. Then
$$\varepsilon_{i}^{-1}({\sf g}_{i}^{2} - {\sf g}_{i}^{1}) = \kappa_{i} +
O(\varepsilon_{i}),$$
where $\kappa_{i} = \varepsilon_{i}^{-1}(\varepsilon_{i}^{2}\kappa_{i}^{2} -
\varepsilon_{i}^{1}\kappa_{i}^{1}) \rightarrow \kappa$, where the convergence is in
$C^{1,\alpha'}$, (in a subsequence). As previously, we need to show that the convergence
is strong, so that $\kappa \neq 0$. This follows from a standard linearization and
bootstrap argument, as preceding \eqref{5.11b}. In more detail, dropping the index $i$,
we have $\Delta_{g^{j}}u^{j} = 0$, so that
$$\Delta_{g^{2}}(u^{1} - u^{2}) = (\Delta_{g^{2}} - \Delta_{g^{1}})u^{1}.$$
In local harmonic coordinates, the right side this equation is on the order of $\varepsilon$
in $C^{1,\alpha}$, and hence by elliptic regularity, $u_{1}-u_{2}$ is on the order of
$\varepsilon$ in $C^{3,\alpha}$. Substituting this in the difference of the static
equations $u^{j}Ric_{g^{j}} = D^{2}_{g^{j}}u^{j}$ and arguing in the same way shows
that the difference $g^{1} - g^{2}$ is then also on the order of $\varepsilon$ in
$C^{3,\alpha}$. This proves the strong convergence.

    It follows that the limit form
\begin{equation}\label{6.8}
\kappa = (h, u'),
\end{equation}
is a non-zero $C^{1,\alpha}$ weak solution of the linearized static vacuum equations
$L(\kappa) = 0$ at $g_{Sch}$ and since $\Pi^{P_{0}}(g_{i}^{1}, u_{i}^{1}) =
\Pi^{P_{0}}(g_{i}^{2}, u_{i}^{2})$, one has
\begin{equation}\label{6.9}
\gamma'_{h} = H'_{h} = 0 \ \ {\rm at} \ \ \partial M,
\end{equation}
where $\gamma'_{h} = h^{T} = h|_{\partial M}$. As discussed following \eqref{5.11b},
elliptic regularity implies that $(h, u')$ is smooth  and so in particular a strong solution.
Below we will use the fact that the data $(h, u')$ are in fact real-analytic up to $\partial M$,
(again by elliptic regularity).

  We claim that
\begin{equation}\label{6.10}
(h, u') = 0 \ \ {\rm on} \ \ M,
\end{equation}
which will give a contradiction. This is of course a linearized version of the
black hole uniqueness theorem. It is possible that \eqref{6.10} can be proved
by linearizing one of the existing proofs of black hole uniqueness in \cite{I},
\cite{R}, \cite{BM}. However, we have not succeeded in doing this and instead
\eqref{6.10} is proved in a manner similar to the proof of Theorem 5.1.

  First, the linearization of \eqref{5.5} gives, at $\partial M$,
$$u' s_{\gamma} = 2\Delta u'.$$
Since $s_{\gamma} > 0$, the maximum principle implies that $u' = 0$ at $\partial M$.
Next we claim $A' = 0$. To see this, the vacuum equations give $uRic = D^{2}u$
and $D^{2}u = N(u)A + (D^{2}u)^{T}$ when evaluated on tangent vectors to
$\partial M$. Taking then the variation and evaluating tangentially gives
$(uRic)' = 0$ so $0 = (D^{2}u)' = (D^{2})'u + D^{2}u'$. The first term
on the right vanishes when evaluating tangentially and hence so does the
second term. This implies $0 = (N(u)A)' = N(u)A' + N(u')A = N(u)A'$. Since
$N(u) = const \neq 0$, it follows that $A' = 0$. Similarly, taking the
variation of the divergence or vector constraint gives $N(u') = const$.

   Clearly $N(u') = m'$, (up to constants). A simple examination of the proof
of black hole uniqueness in \cite{R} applied to an Einstein deformation as in
\eqref{6.5} and satisfying \eqref{6.9}, shows easily that $N(u') = 0$ at $\partial M$.
(One does not obtain any further information, since the bulk data in the Robinson proof,
via divergence identities, are quadratic in the deviation from Schwarzschild).

  Thus the variations $(\gamma', u', A', N(u'))$ of all the Cauchy data are trivial.
As in the proof of Theorem 4.1, we use a bootstrap argument to prove that
the data $(h, u')$ vanish to infinite order at $\partial M$, in geodesic gauge.

  Thus, using geodesic normal coordinates near $\partial M$, write
$$g = dt^{2} + g_{t},$$
where $t(x) = dist(x, \partial M)$. We may assume, (by adding an infinitesimal
deformation of the form $\delta^{*}V$ if necessary), that $h$ preserves this gauge,
so that $h_{0\alpha} = 0$, i.e. $h(N, \cdot) = 0$, $N = \partial_{t}$, near
$\partial M$. By the
discussion above, we have $u' = N(u') = 0$ at $\partial M$ and similarly $h =
\nabla_{N}h = 0$ at $\partial M$, so that
\begin{equation}\label{6.11}
u' = O(t^{2}) \ \ {\rm and} \ \  h = O(t^{2}).
\end{equation}

  The variation of the potential equation $\Delta_{M}u = 0$ gives
\begin{equation}\label{6.12}
\Delta u' = -\Delta'u = \langle D^{2}u, h \rangle - \langle \beta (h), du
\rangle,
\end{equation}
where $\beta$ is the Bianchi operator, (cf.~\cite{Be}). Since $\beta (h) = 0$
at $\partial M$, this gives $\Delta u' = 0$ at $\partial M$ and hence $NN(u') = 0$
at $\partial M$, so that
\begin{equation}\label{6.13}
u' = O(t^{3}).
\end{equation}
Next the linearization of the Riccati equation gives $(\nabla_{N}A)' = -(A^{2})' -
(R_{N})' = -(R_{N})'$ at $\partial M$. One computes $*R_{N} = -(Ric^{T}) =
- u^{-1}(D^{2}u)^{T}$, so $(*R_{N})' =  u^{-2}u'D^{2}u - u^{-1}(D^{2})'u -
u^{-1}D^{2}u'$, as a form on $T(\partial M)$. It follows from \eqref{6.13}
that $(*R_{N})' = 0$ at $\partial N$ and hence $(\nabla_{N}A)' = 0$ so that
\begin{equation}\label{6.14}
h = O(t^{3}).
\end{equation}
Next taking the normal derivative of \eqref{6.12} gives
$$\Delta N(u') = \langle \nabla_{N}D^{2}u, h \rangle + \langle D^{2}u,
\nabla_{N}h \rangle +  \langle \nabla_{N}\beta (h), du \rangle +
\langle \beta (h), \nabla_{N}du \rangle,$$
which vanishes at $\partial M$ and hence $u' = O(t^{4})$. Substituting this in
the linearized Riccati equation above and using previous estimates gives
$h = O(t^{4})$, and so on. It follows that $(h, u')$ vanish to infinite order
at $\partial M$. Since $(h, u')$ are real-analytic up to $\partial M$, (in
geodesic gauge), this implies that $h = u' = 0$ on $M$, i.e.~\eqref{6.10} holds.

   It remains only to prove that there are regular values $(\gamma, H)$ of $\Pi^{P_{0}}$ near
$(\gamma_{+1}, 0)$ whose inverse image under $\Pi^{P_{0}}$ is non-empty. However,
\eqref{6.10} proves that $Ker D\Pi_{B} = 0$ at the Schwarzschild metric with boundary
data $(\gamma_{+1}, 0)$. By continuity and the black hole uniqueness theorem, it
follows that $Ker D\Pi_{B} = 0$, for all solutions with boundary data $(\gamma, H)$
near $(\gamma_{+1}, 0)$ in $Im \Pi^{P_{0}}$. Thus, all such boundary data are regular
values of $\Pi^{P_{0}}$. This completes the proof.

{\endproof}

\begin{remark}\label{r6.3}
{\rm  Although the proof of Theorem 5.2 implies that $D\Pi_{B}$ has trivial kernel at
the Schwarzschild metric $g_{Sch}$, one does not expect this to be the case for the
cokernel. In fact, one expects that $Coker D\Pi_{B}$ is infinite dimensional, in that
any boundary variation of the form $(k, 0)$, where $k$ is a variation of the boundary
metric, is not tangent to a curve of static metrics with $H = 0$ at $\partial M$.
This amounts to the linearized version of the no-hair theorem, (which has not been
proved as far as we are aware). In particular, we expect $D\Pi_{B}$ is not
Fredholm at $g_{Sch}$. }
\end{remark}

\begin{remark}\label{r6.4}
{\rm The proof of Theorem 5.2 above shows that $Im \Pi^{P_{0}}$ and hence $Im \Pi_{B}$
contains regular values. In particular, this means that the image $Im \Pi_{B}$ has non-empty
interior and is in fact an open map near the Schwarzschild metric with boundary the
horizon $r = 2m$. We point out that it has been an open question as to whether
$\Pi_{B}$ has {\it any} regular values. For instance, the analysis of Miao in \cite{M1}
shows that the exterior of the standard round unit ball $B^{3}(1)$ in flat
${\mathbb R}^{3}$ is a critical point of $\Pi_{B}$. In fact it is still unknown whether
$Im \Pi_{B}$ has open interior near such standard boundary data
$(\gamma_{+1}, 2)$. }
\end{remark}

  Remark 5.3 indicates how little has been understood, and how much remains to be
understood, regarding the global behavior of the boundary map $\Pi_{B}$. The proof
of the almost properness of $\Pi^{o}$ in Theorem 1.2 uses the outer-minimizing
property \eqref{1.10} in two key, but essentially independent ways. First, it is used
to obtain the a priori curvature and related estimates discussed in Section 4,
i.e.~to prevent blow-up behavior at the boundary. For this, only a small-scale
version of \eqref{1.10} is necessary, in that one needs stability of the area of the
boundary only to 4th order at $\partial M$, and this only in small discs in $\partial M$.
Second, it is used to keep the boundary $\partial M$ properly embedded in $M$,
i.e.~to prevent the passage from embedded to immersed behavior. Again, as
mentioned in Remark 4.4, only a local version of \eqref{1.10} is needed for this.

   For further progress, an important issue is to find some condition on the boundary
data $(\gamma, H)$ which ensures that the two properties above hold. A natural
question is whether positivity of the Gauss and mean curvatures is sufficient.

\bibliographystyle{plain}

\end{document}